\documentclass[10pt,reqno]{article}

\usepackage{latexsym, amsmath, amssymb}
\usepackage{graphicx,epsfig}

\usepackage{amsfonts}
\usepackage{amscd}
\usepackage{amsbsy}
\usepackage{bm}
\usepackage{color}
\usepackage{cite}

\newtheorem{theorem}{Theorem}[section]

\newtheorem{lemma}{Lemma}[section]

\newtheorem{proposition}{Proposition}[section]

\numberwithin{equation}{section}

\newcommand{\Ga}{\alpha}
\newcommand{\Gb}{\beta}

\newcommand{\Gg}{\gamma}
\newcommand{\Gs}{\sigma}
\newcommand{\Ggs}{s}
\newcommand{\Gl}{\lambda}
\newcommand{\GD}{\Delta}

\newcommand{\GG}{\Gamma}
\newcommand{\GL}{\Lambda}

\newcommand{\Bx}{{\bf x}}
\newcommand{\By}{{\bf y}}
\newcommand{\Bz}{{\bf z}}
\newcommand{\BM}{{\bf M}}
\newcommand{\BV}{{\bf V}}

\newcommand{\qed}{\ \ensuremath{\square}}
\newcommand{\ds}{\displaystyle}
\newcommand{\pf}{\medskip \noindent {\sl Proof}. ~ }
\newcommand{\p}{\partial}
\renewcommand{\a}{\alpha}

\newcommand{\pd}[2]{\frac {\p #1}{\p #2}}

\newcommand{\eqnref}[1]{(\ref {#1})}
\newcommand{\na}{\nabla}
\newcommand{\Om}{\Omega}
\newcommand{\ep}{\epsilon}

\newcommand{\RR}{\mathbb{R}}

\newcommand{\NN}{\mathbb{N}}

\newcommand{\la}{\langle}
\newcommand{\ra}{\rangle}

\newcommand{\Scal}{\mathcal{S}}

\newcommand{\Kcal}{\mathcal{K}}

\newcommand{\beq}{\begin{equation}}
\newcommand{\eeq}{\end{equation}}

\def\nm{\noalign{\medskip}}

\begin{document}
\title{Numerical implementation for reconstruction of inhomogeneous conductivities via Generalized Polarization Tensors
\thanks{Fang is supported by NSF grants  No. 70921001, No. 71210003 and No. 71271219. Deng is supported by NSF grants No. NSFC11301040.}}
\author{Xiaoping Fang \thanks{Postdoctoral, Management Science and Engineering Postdoctoral Mobile Station, School of Business; School of Mathematics and Statistics, Central South University, Changsha, Hunan 410083, P. R. China. Email: fxpmath@csu.edu.cn} \quad Youjun Deng \thanks{Corresponding author. School of Mathematics and Statistics, Central South University, Changsha, Hunan 410083, P. R. China. Email: youjundeng@csu.edu.cn, dengyijun\_001@163.com}
}

\date{}
\maketitle

\begin{abstract}
This paper deals with numerical methods for reconstruction of inhomogeneous conductivities. We use the concept of Generalized Polarization Tensors, which were introduced in \cite{ADKL12}, to do reconstruction. Basic resolution and stability analysis are presented. Least square norm methods with respect to Generalized Polarization Tensors are used for reconstruction of conductivities. Finally, reconstruction of three different types of conductivities in the plane is demonstrated.
\end{abstract}

\noindent {\footnotesize {\bf Mathematics subject classification
(MSC2000):} 35R30, 35C20}


\section{Introduction}
It is known that  Generalized Polarization Tensors (GPTs) carry geometric information about the inclusion. They uniquely determine the conductivity distribution \cite{ADKL12,AK02}.
More importantly, in some sense, GPTs determine the conductivity hierarchically from lower order to higher order. To explain this,
we consider the homogeneous conductivity problem in $\RR^d$, $d=2,3$:
\begin{equation}\label{condeqn}
\left\{
\begin{array}{ll}
\nabla \cdot (1 + (c-1) \chi(\Om)) \nabla u =0 & \mbox{in } \RR^d, \\
\nm u(\Bx)-h(\Bx)=O(|\Bx|^{1-d}) & \mbox{as } |\Bx| \rightarrow \infty,
\end{array}
\right.
\end{equation}
where $\Om$ is the inclusion embedded in $\RR^d$ with a Lipschitz boundary, $\chi(\Om)$ is the characteristic function of $\Om$.
The positive constant $c>0$ is the conductivity of the inclusion. $h$ is a harmonic
function in $\RR^d$ representing the background electrical potential, and the solution $u$ to the problem represents the
perturbed electrical potential. The perturbation $u-h$ due to the presence of the conductivity inclusion $\Omega$ admits the
following asymptotic expansion as $|\Bx| \to \infty$ (see \cite{ADKL12,book2}):
 \beq\label{uhasymphom}
 u(\Bx)-h(\Bx) = \sum_{|\Ga|, |\Gb| \ge 1} \frac{(-1)^{|\Gb|}}{\Ga ! \Gb !}
 \p^\Ga h(0) M_{\Ga\Gb}(c, \Om) \p^\Gb \Gamma(\Bx),
 \eeq
where $\GG$ is the fundamental solution of the Laplacian which has the form
\beq
\Gamma(\Bx)=\left\{
\begin{array}{cl}
 \frac{1}{2\pi}\ln|\Bx|, & d=2,\\
 -\frac{1}{4\pi}\frac{1}{|\Bx|}, & d=3,
\end{array}
\right.
\eeq
and $\a$ and $\beta$ are the multi-indexes.
The building blocks  $M_{\Ga\Gb}(c,\Om)$ for the  asymptotic expansion  \eqnref{uhasymphom} are called the GPTs.
The leading order GPT (called the polarization tensor (PT)), $\{ M_{\Ga\Gb}(c,\Om):  |\Ga|, |\Gb| = 1 \}$, provides the equivalent ellipse
(ellipsoid) which represents overall property of the inclusion \cite{AKKL03, BHV03}. The concept of GPTs was first introduced
by Ammari et al \cite{book2} for electrical impedance imaging of small conductivity inclusions, then extended to various areas such as
elastic imaging \cite{AKLL13}, the theory of dilute composite materials, electric magnetic wave imaging and so on.
The GPTs carry geometric information about the inclusion. For example, the whole set of GPTs,
$\{ M_{\Ga\Gb}(c,\Om): \Ga , \Gb \geq 1 \}$, determines the conductivity $c$ and the inclusion $\Om$ uniquely \cite{AK02}.
Moreover, there are important analytical and numerical
studies which show that finer details of the shape can be recovered using higher-order GPTs \cite{AGKLY11, book2, handbook, AKLZ12}.
The GPTs even carry topology information of the inclusion \cite{AGKLY11}, although accurate topological information carried
by different order of GPTs are yet to be investigated. The contracted GPTs (CGPTs), which are harmonic combinations of GPTs
\cite{AKLL11}, turn out to be more efficient in reconstruction of conductivity, since they associates with different
frequency components of the boundary measurements.

This paper is concerned with the numerical implementation of the reconstruction of inhomogeneous conductivities in the plane,
which means that $c$, in \eqnref{uhasymphom},  is replaced by some non-constant function $\Gs$.
The famous Calder\'on problem \cite{Cald80} is to determine $\Gs$ from the knowledge of DtN (Dirichlet-to-Neumann) map .
We only consider the isotropic conductivity, as it is proved that all the boundary measurements (the Dirichlet-to-Neumann map) could not
uniquely determine the anisotropic conductivity by using change of coordinates \cite{ALPa05}. For unique determination
of $\Gs$ from the DtN map, we refer the classical results in three and higher
dimensional smooth conductivities by Sylvester and Uhlmann \cite{SyUh87}. In two dimensional problem, uniqueness was established
for piecewise analytic conductivities by Kohn and Vogelius \cite{KoVo84,KoVo85} and for generic by Sun and Uhlmann \cite{SuUh91}.
Later, A. Nachman \cite{Nach96} proved the uniqueness for conductivities with two derivatives. The uniqueness for two dimensional
bounded measurable conductivity was finally solved by Astala and P$\ddot{a}$iv$\ddot{a}$rinta in 2006 by using the quasi-conformal
mappings \cite{AP06}.
For reconstruction, we use the notion of GPTs for inhomogeneous conductivities in \cite{ADKL12}, where
some special numerical example is shown, i.e., the
reconstruction of radially symmetric conductivity.
In this paper, we consider the reconstruction of the conductivities of more general cases. To do this, we define an operator
which contains all the information of GPTs (and contracted GPTs). The operator is related to the
Neumann-to-Dirichlet (NtD) map and we show that the stability result of reconstruction of inhomogeneous
conductivities by using this operator can be obtained similarly to that by using the NtD map. We also show that
the eigenvalues of the operator is actually the linear combinations of contracted GPTs in some special cases.
The basis for reconstruction of the inhomogeneous conductivity is to minimize the discrepancy functional between the
reconstructed GPTs from observations and the GPTs related to undetermined conductivity. It is shown in \cite{ABGKW13}
 that in the full-view case, the reconstruction problem of GPTs from boundary data
has the remarkable property that low order GPTs are not affected by the error caused by the instability of higher-orders
in the presence of measurement noise.

The organization of this paper is as follows. In section 2, we present the contracted GPTs for the inhomogeneous conductivity distribution which has been
defined in \cite{ADKL12}. We also show the far field expansion of the perturbed potential by using contracted GPTs. We then define an operator
in section 3 and prove some properties of the operator. The stability result is also obtained by using this operator to reconstruct the conductivity
$\Gs$. In section 4, we present the optimization method for reconstruction of the GPTs by using Multi Static Response(MSR) matrix. Section 5 deals with the reconstruction of the conductivity distribution by using optimization methods. Numerical experiments for reconstructing three different types of conductivities are presented in section 6.

\section{Contracted GPTs}
Let $\Gs$ be a bounded measurable function in $\RR^d$, $d=2,3$ such that
 \beq\label{glonetwo}
 c^{-1} \le \Gs \le c
 \eeq
for positive constant $c>1$. In this paper, we only consider the reconstruction of
H$\ddot{o}$lder continuous conductivity $\Gs$, i.e.,
$\Gs_i\in C^s(D)$ with $s>0$. But we point out that most of the theories presented are
suitable for any $L^{\infty}$ conductivity unless otherwise pointed out(see \cite{ADKL12}).
Let $h$ be the harmonic function in $\RR^d$. Suppose $h$ is perturbed in the environment with conductivity
$\Gs$ and $u$ is the perturbed potential, then $u$ is the solution to
\begin{equation}\label{eq:101}
\left\{
\begin{array}{ll}
\nabla \cdot \sigma \nabla u =0 & \mbox{in } \RR^d, \\
\nm u(\Bx)-h(\Bx)=O(|\Bx|^{1-d}) & \mbox{as } |\Bx| \rightarrow \infty.
\end{array}
\right.
\end{equation}
Let $D$ be a bounded domain in $\RR^d$ with a
$\mathcal{C}^{1,\eta}$-boundary $\partial D$ for some $0<\eta<1$.
We assume that $D$ is such that
 \beq\label{suppD}
 \mbox{supp}\, (\Gs-1) \subset D.
 \eeq
 Let $H^{s}(\p D)$, for $s\in\RR$,  be the usual $L^2$-Sobolev space and let
$$H^{s}_0 (\p D):=\left\{\phi \in H^{s}(\p D) \Big| \int_{\p D} \phi =0\right\}.$$
For $s=0$, we
use the notation $L^2_0(\partial D)$. We define the Neumann-to-Dirichlet (NtD) map $\GL_{\Gs} : H^{-1/2}_0(\p D) \to H^{1/2}_0(\p D)$ as
 \beq
 \GL_{\Gs}[g]:= u|_{\p D},
 \eeq
where $u$ is the solution to
\beq\label{sigmaeqn}
\left\{
\begin{array}{ll}
\nabla \cdot \Gs\nabla u=0 & \mbox{in }  D, \\
\ds \Gs\frac{\partial u}{\partial \bm\nu} = g & \mbox{on } \ds \partial D \quad \left(\int_{\p D} u=0\right)
\end{array}
\right.
\eeq
for $g \in H^{-1/2}_0(\p D)$. We mention that lots of references use the Dirichlet-to-Neumann (DtN) map in stead
of NtD map for analysis (see, e.g., \cite{ChYa04,CLNa05,KSVW08}). Since we impose the operator $\GL_{\Gs}$ on
$H^{-1/2}_0(\p D)$, the inverse operator $\GL_{\Gs}^{-1}$ is the DtN map acting from $H^{1/2}_0(\p D)$ to
$H^{-1/2}_0(\p D)$.
The operator $\GL_1$ is the NtD map when $\Gs\equiv 1$. In the following
we define the operator $I_0$ which maps from $W(\p D)$ to $W_0(\p D)$ by
$$
I_0 \varphi := \varphi - \int_{\p D} \varphi(\By) ds(\By), \quad \varphi \in W(\p D)
$$ where $W$ can be any Hilbert space.
Let $\Lambda^e$ be the NtD map for the exterior problem:
$$ \Lambda^e[g]:=I_0 u|_{\p D}$$
where $u$ is the solution to
\begin{equation}\label{eq:extop}
\left\{
\begin{array}{ll}
\GD u =0 & \mbox{in }  \RR^d \setminus \overline{D}, \\
\ds \frac{\p u}{\p \bm\nu} \Big|_+=g & \mbox{on } \p D,\\
u(\Bx) =O(|\Bx|^{1-d})  & \mbox{as } |\Bx| \rightarrow \infty.
\end{array}
\right.
\end{equation}
We shall restrict the operator $\GL^e$ on $H_0^{-1/2}$. Then we have
$$
\GL^e[g]=I_0\Scal_D(\frac{1}{2}I+\Kcal_D^*)^{-1}[g]
$$
for any $g\in H^{-1/2}_0(\p D)$, where $\Scal_D$ is the single layer potential defined by
 \beq\label{eq:sglyer}
 \Scal_D[\phi](\Bx)=\int_{\p D} \GG(\Bx-\By) \phi(\By)
 d\Ggs(\By), \quad x \in \RR^d.
 \eeq
 and $\Kcal_D^*$ is the adjoint operator of Poincar\'e-Neumann operator $\Kcal_D$ defined by
 $$
\Kcal_D[\phi](\Bx)= \int_{\p D} \frac{\partial \GG}{\partial
\bm\nu_y}(\Bx-\By) \phi(\By) d\Ggs(\By).
$$
Thus there holds $\GL^e=I_0\Scal_D(1/2 I +\Kcal_D^*)^{-1}$ from $H_0^{-1/2}$ to $H_0^{1/2}$.
For $d=2$, the contracted generalized polarization tensors (CGPTs) are defined as follows (see
\cite{AKLL11}):
 \begin{align}
 \ds M_{mn}^{cc}= M_{mn}^{cc}(\Gs,D):= \int_{\p D} r_y^m\cos m\theta_y \, \GL_1^{-1} (\GL_1- \GL_{\Gs})
 [g_n^c] (\By) \, d\Ggs(\By), \label{defmcc}\\
 \ds M_{mn}^{cs}= M_{mn}^{cs}(\Gs,D):=\int_{\p D} r_y^m\cos m\theta_y \, \GL_1^{-1} (\GL_1- \GL_{\Gs}) [g_n^s] (\By) \, d\Ggs(\By), \\
 \ds M_{mn}^{sc}= M_{mn}^{sc}(\Gs,D):= \int_{\p D} r_y^m\sin m\theta_y \, \GL_1^{-1} (\GL_1- \GL_{\Gs}) [g_n^c] (\By) \, d\Ggs(\By), \\
 \ds M_{mn}^{ss}= M_{mn}^{ss}(\Gs,D):= \int_{\p D} r_y^m\sin m\theta_y \, \GL_1^{-1} (\GL_1- \GL_{\Gs}) [g_n^s] (\By) \, d\Ggs(\By),
 \end{align}
with $\Bx=(r\cos\theta, r\sin\theta)$ and $g_n^s$ and $g_n^c$ satisfy
\beq
\label{eq:gn}
g_n^s= (\GL_{\Gs}-\GL^e)^{-1}(\GL_1-\GL^e)[\na (r^{n}\sin n\theta) \cdot \bm\nu],  \, \, g_n^c= (\GL_{\Gs}-\GL^e)^{-1}(\GL_1-\GL^e)[\na (r^{n}\cos n\theta)\cdot \bm\nu],
\eeq
where $\bm\nu$ is unit outward normal to the domain $D$.
We have the far field expansion of the perturbation $u-h$ to \eqnref{eq:101} in the form of contracted GPTs.
\begin{theorem}
\label{th:expan2} Let $u$ be the solution to \eqnref{eq:101} with
$d=2$. If $h$ admits the expansion
\beq\label{Hexp}
 h(\Bx) = h(0)+\sum_{n=1}^\infty r^n \bigr(a_n^c\cos n\theta + a_n^s\sin n\theta \bigr)
 \eeq
then we have
\begin{align}
(u-h)(\Bx) & = -\sum_{m=1}^\infty\frac{\cos m\theta}{2\pi mr^m}\sum_{n=1}^\infty
\bigr(M_{mn}^{cc}a_n^c + M_{mn}^{cs}a_n^s \bigr) \nonumber \\
& \qquad -\sum_{m=1}^\infty\frac{\sin m\theta}{2\pi m
r^m}\sum_{n=1}^\infty \bigr(M_{mn}^{sc}a_n^c + M_{mn}^{ss}a_n^s
\bigr), \label{expan2}
\end{align}
which holds uniformly as $|\Bx| \to \infty$.
\end{theorem}

For $d=3$, we also have
\begin{theorem}
\label{th:expan3} Let $u$ be the solution to \eqnref{eq:101} with
$d=3$. If $h$ admits the expansion
\begin{align}
h(\Bx)=h(0) + \sum_{m=1}^\infty \sum_{k=-m}^m a_{mk} r^m Y_m^k(\theta, \varphi), \label{hdecom}
\end{align}
then we have
\beq\label{expan3} (u-h)(\Bx) =- \sum_{m=1}^\infty
\sum_{k=-m}^{m}\sum_{n=1}^\infty \sum_{\ell=-n}^n
\frac{a_{mk}M_{mkn\ell}}{(2n+1)r^{n+1}}Y_n^\ell(\theta,
\varphi)\quad\mbox{\rm as } |x| \to \infty, \eeq where the GPT
$M_{mnk\ell}=M_{mnk\ell}[\Gs]$ is defined by \beq
M_{mkn\ell}:=\int_{\partial B} Y_m^k (\theta', \varphi') r'^m
\GL_1^{-1} (\GL_1- \GL_{\Gs})[g_{nl}] (r',\theta', \varphi') \,
ds \eeq
where $g_{nl}:=\Gs\frac{\p u_{nl}}{\p \bm\nu}\Big|_{-}$ and $u_{nl}$ is the solution to \eqnref{eq:101} with $h$
replaced by $r_{n} Y_n^l(\theta,\varphi)$.
\end{theorem}

For the sake of simplicity, we denote by $M_{mn}$ the contracted GPTs while omitting the $c$ and $s$ for the superscripts
in two dimensional problem and subscripts $k$ and $l$ in three dimensional problem.
In the sequel, we denote by $h_m$ the harmonic function with degree $m$ in $\RR^d$. $h_m$ is the linear combination of
the harmonic functions of order $m$. Then $M_{mn}$ can be written as
\beq\label{eq:conGPTs}
M_{mn}:=M_{mn}(\Gs)= \int_{\p D} h_m(\By)  \phi_n (\By) \, d\Ggs(\By),
\eeq
where
\beq\label{eq:phin}
\phi_n=\GL_1^{-1} (\GL_1 - \GL_\Gs) (\GL_\Gs - \GL^e )^{-1}
  (\GL_1 - \GL^e ) \left[ \pd{h_n}{\bm\nu} \Big|_{\p D} \right].
\eeq
For the invertibility of operator $\GL_\Gs - \GL^e $ see Lemma \ref{le:bnd}.

\section{Stability analysis}
In this section, we shall only consider the two dimensional problem while the same strategy can
be used for the analysis of the three dimensional problem.
To simplify the stability analysis we define the operator $M_{\Gs}: H^{-1/2}_0(\p D) \to H^{-1/2}_0(\p D)$ by
\beq\label{eq:GPTop}
M_{\Gs}:= \GL_1^{-1}(\GL_1 - \GL_\Gs) (\GL_\Gs - \GL^e )^{-1}(\GL_1 - \GL^e )
\eeq
The operator $M_{\Gs}$ can actually be treated as the GPTs operator, i.e., we have the following relation
\beq
M_{mn}=\la h_m, M{\Gs}\left[\frac{\p h_n}{\p \bm\nu}\right]\ra_{H^{1/2},H^{-1/2}}=:\la h_m , \frac{\p h_n}{\p \bm\nu}\ra_{M_{\Gs}}
\eeq
where $\la \cdot ,\cdot \ra_{H^1/2,H^{-1/2}}$ is the bilinear product. If $\Gs$ is a constant in $D$ then
it is not difficult to verify that
$$
M_{\Gs}=\left(\frac{\Gs+1}{2(\Gs-1)}I-\Kcal_D^*\right)^{-1}.
$$
In the sequel, we suppose $C$ is a common
positive constant which may change in each occurrence.
We give some properties for $M_{\Gs}$, but before this we present some primary lemmas.
\begin{lemma}[Lemma 5.2 in \cite{ADKL12}]
\label{le:est_op}
There is a constant $C$ such that
\beq\label{DtNstability}
\|\GL_{\Gs_1}-\GL_{\Gs_2}\|  \leq  C \|\Gs_1-\Gs_2\|_{L^{\infty}(D)} .
\eeq
\end{lemma}
\begin{lemma}
\label{le:est4}
Let $u_j$, $j=1,2$ be the solution of the following problem
$$
\left\{
\begin{array}{ll}
\na \cdot \Gs_j \na u_j=0 & \mbox{in} \ \ \RR^2\\
u_j -h = O(|\Bx|^{-1}) & \mbox{as} \ \ |\Bx|\rightarrow \infty,
\end{array}
\right.
$$
where $\Gs_j$ satisfy \eqnref{glonetwo} and $supp(\Gs_j-1) \subset D$, $j=1,2$,
then there holds
$$
\|\na(u_1-u_2)\|^2_{L^2(\RR^2\setminus D)} \leq C \|\na h\|^2_{L^2(D)}\|\Gs_1-\Gs_2\|_{L^{\infty}(\RR^2)}.
$$
\end{lemma}
\pf By using the representation of the solutions in \cite{ADKL12} we know that $u_1$ and $u_2$ has the following form in $\RR^2\setminus D$
 $$
 u_j= h+ \Scal_D[\phi_j], \quad \phi_j=M_{\Gs}[\frac{\p h}{\p \bm\nu}], j=1,2
 $$
Furthermore, by Lemma \ref{le:est_op} one can get
\begin{eqnarray*}
& &\int_{\RR^2\setminus D}|\na(u_1-u_2)|^2 d\Bx  =  -\int_{\p D} \frac{\p(u_1-u_2)}{\p \bm\nu}\Big|_+ (u_1-u_2) ds  \\
& = & \int_{\p D} (\phi_1-\phi_2)\GL_1 (\GL_1-\GL^e)^{-1}\GL_1 (\GL_1-\GL^e)^{-1}\GL_1[\phi_1-\phi_2] ds \\
& = & \Big\|\left((\GL_1+\GL_{\Gs_1})^{-1}(\GL_1-\GL_{\Gs_1}) -(\GL_1+\GL_{\Gs_2})^{-1}(\GL_1-\GL_{\Gs_2})\right)\Big[\frac{\p h}{\p \bm\nu}\Big]\Big\|^2_{H^{-1/2}(\p D)} \\
& \leq & 2\Big\|\left((\GL_1+\GL_{\Gs_1})^{-1}(\GL_{\Gs_2}-\GL_{\Gs_1})\right)\Big[\frac{\p h}{\p \bm\nu}\Big] \Big\|^2_{H^{-1/2}(\p D)} \\
& & + 2\Big\| \left(\left((\GL_1+\GL_{\Gs_1})^{-1}-(\GL_1+\GL_{\Gs_2})^{-1}\right)(\GL_1-\GL_{\Gs_2})\right)\Big[\frac{\p h}{\p \bm\nu}\Big] \Big\|^2_{H^{-1/2}(\p D)} \\
& \leq & C\|\na h\|^2_{L^2(B)} \|\Gs_2-\Gs_1\|_{L^{\infty}(D)}
\end{eqnarray*}
which completes the proof.
\qed

\begin{lemma}\label{le:radsol}
Suppose $\Gs$ is radially symmetric. Let $u_m$ be the solution of the following problem
$$
\left\{
\begin{array}{ll}
\na \cdot \Gs \na u_m=0 & \mbox{in} \ \ \RR^2\\
u_m -r^m\cos m\theta \rightarrow 0 & \mbox{as} \ \ r\rightarrow \infty,
\end{array}
\right.
$$
where $m \in \NN$ and $D$ is defined as usual then the solution
is unique and has the form $u_m=(r^m+\frac{b_m}{r^m})\cos m\theta$ for $\Bx\in \RR^2\setminus D$.
\end{lemma}
\pf Since $\Gs$ is radially symmetric, for any $\Gs$ we can find piecewise constant conductivities $\{\Gs_i\}$ such that
$\|\Gs_i-\Gs\|_{L^{\infty}(\RR^2)}\rightarrow 0$. It is easy to see that the conductivity problem
$$
\left\{
\begin{array}{ll}
\na \cdot \Gs_i \na u_m^{(i)}=0 & \mbox{in} \ \ \RR^2\\
u_m^{(i)} -r^m\cos m\theta \rightarrow 0 & \mbox{as} \ \ r\rightarrow \infty,
\end{array}
\right.
$$
has a unique solution with $u_m^{(i)}=(r^m+\frac{b_m^{(i)}}{r^m})\cos m\theta$ for $\Bx\in \RR^2\setminus D$. On the
other hand, by Lemma \ref{le:est4}, we have
$$
\|\na(u_m^{(i)}-u_m)\|^2_{L^2(\RR^2\setminus D)} \leq C \|\na (r^m \cos m\theta)\|^2_{L^2(B)}\|\Gs_i-\Gs\|_{L^{\infty}(\RR^2)}.
$$
Thus
$$
\lim_{i\rightarrow \infty} \|\na(u_m^{(i)}-u_m)\|^2_{L^2(\RR^2\setminus D)}=0.
$$
So $u_m$ must have the form $(r^m+\frac{b_m}{r^m})\cos m\theta+c$ for $\Bx\in \RR^2\setminus D$. Since
$u_m-r^m\cos m\theta \rightarrow 0$ we immediately get $c=0$ and the proof is complete.
\qed

\begin{proposition}\label{prop:prop}
There holds the following for the operator $M_{\Gs}: H^{-1/2}_0(\p D) \to H^{-1/2}_0(\p D)$
\begin{enumerate}
  \item[(i)] $M_{\Gs}$ is self-adjoint.
  \item[(ii)]
If $\Gs$ is radially symmetric, and $D$ is a unit disk, then $M_{\Gs}$ has eigenfunctions $\sqrt{\frac{m}{\pi}}\cos m\theta$ and
corresponding eigenvalues $M_{mm}(\Gs)$ for all $m\in \NN$, i.e.,
\beq\label{eq:eigenM}
M_{\Gs}\Big[\sqrt{\frac{m}{\pi}}\cos m\theta\Big]= M_{mm}(\Gs) (\frac{\cos m\theta}{\sqrt{m\pi}}).
\eeq
\end{enumerate}
\end{proposition}
\pf (i). It follows from
$$
(\GL_1 - \GL_\Gs) (\GL_\Gs - \GL^e )^{-1}(\GL_1 - \GL^e )= -(\GL_1 - \GL^e) + (\GL_1 - \GL^e) (\GL_\Gs - \GL^e )^{-1}(\GL_1 - \GL^e)
$$
and the fact that the NtD map is self-adjoint that $M_{\Gs}$ is self-adjoint.

(iii). It is easy to verify that $\sqrt{\frac{m}{\pi}}\cos m\theta$ and $\frac{\cos m\theta}{\sqrt{m\pi}}$, $m\in \NN$ are the orthogonal basis in $H_0^{-1/2}(\p D)$ and $H_0^{1/2}(\p D)$ when $D$ is a disk, respectively
. Firstly we prove that $M_{mn}=0$ for $m\neq n$. By Lemma \ref{le:radsol} we suppose the solution to
$$
\left\{
\begin{array}{ll}
\na \cdot \Gs \na u_m=0 & \mbox{in} \ \ \RR^2\\
u_m -r^m\cos m\theta \rightarrow 0 & \mbox{as} \ \ r\rightarrow \infty,
\end{array}
\right.
$$
has the form $u_m= (r^m+ \frac{b_m}{r^m})\cos m\theta$. Thus we get
\beq
\label{eq:NtDmap}
\GL_{\Gs}[\cos m\theta]= \frac{1}{m}\frac{1+b_m}{1-b_m}\cos m\theta.
\eeq
Since $D$ is a disk, $\GL^e=-\GL_1$. It then follows that
$$
M_{mn}= 2\int_{\p B} \sqrt{\frac{m}{\pi}}\cos m\theta(\GL_{1}-\GL_{\Gs}) (\GL_1+\GL_{\Gs})^{-1}\GL_1\Big[\sqrt{\frac{m}{\pi}}\cos m\theta\Big]ds
= -2 b_m \delta_{mn}
$$
and \eqnref{eq:eigenM} follows immediately.
\qed

Let $\Gs_1$ and $\Gs_2$ be two different conductivity distribution, we have the following relation
\beq\label{eq:reM}
M_{\Gs_1}-M_{\Gs_2}= \GL_1^{-1}(\GL_1-\GL^e)(\GL_{\Gs_1}-\GL^e)^{-1}(\GL_{\Gs_2}-\GL_{\Gs_1})(\GL_{\Gs_2}-\GL^e)^{-1}(\GL_1-\GL^e).
\eeq
With \eqnref{eq:reM} on hand, we have the following result
\begin{lemma}\label{le:bnd}
\beq\label{eq:Gsbound}
C_2 \|M_{\Gs_1}-M_{\Gs_2}\|\leq \|\GL_{\Gs_1}-\GL_{\Gs_2}\|\leq C_1 \|M_{\Gs_1}-M_{\Gs_2}\|,
\eeq
where $C_1$ and $C_2$ are constants depend on $\Gs_i$, $i=1,2$.
\end{lemma}
\pf  We first prove that $\GL_{\Gs}-\GL^e$ is bounded invertible for any $\Gs$ satisfying \eqnref{glonetwo}.
For any $g\in H_0^{-1/2}(\p D)$, let $u_1$ and $u_2$ be the solution to \eqnref{sigmaeqn} and \eqnref{eq:extop}, respectively.
Then we have
\begin{align*}
\int_{\p D}(\GL_{\Gs}-\GL^e)[g]g ds & = \int_{\p D}\GL_{\Gs}[g]gds -\int_{\p D}\GL^e[g]gds \\
& = \int_D \Gs |\nabla u_1|^2 ds - \int_{\p D}u_2\frac{\p u_2}{\p \bm\nu}ds \\
& = \int_D \Gs |\nabla u_1|^2 d\Bx + \int_{\RR^d\setminus \overline{D}}|\nabla u_2|^2 d\Bx
\end{align*}
Let
$h$ be the harmonic function in $D$ with $\pd{h}{\nu} = g$ on $\p D$. Then we have the estimate
$$
\int_{D} |\na h|^2d\Bx=\int_{D} \Gs \na u_1\cdot \nabla h d\Bx \leq \frac{\epsilon}{2} \int_{D} |\na h|^2 d\Bx+ \frac{1}{2\epsilon}\int_{D} \Gs^2|\na u_1|^2 d\Bx,
$$
thus by choosing appropriately a small $\epsilon$ we have
$$
\int_D \Gs |\nabla u_1|^2 d\Bx \geq c^{-1} \int_D \Gs^2 |\nabla u_1|^2 d\Bx \geq C\int_{D} |\na h|^2 d\Bx=C\| g\|_{H^{-1/2}(\p D)}^2.
$$
We have shown that $\GL_\Gs$ is strictly positive definite operator and $\GL^e$ is negative definite operator. The injectivity
of the operator $\GL_\Gs-\GL^e$ is then obtained by
$$\|(\GL_\Gs-\GL^e)[g]\|\geq C\|g\|_{H^{-1/2}}(\p D).$$
To prove that The $\GL_\Gs-\GL^e$ is surjective is equivalent to prove the existence of the following equation
\begin{equation}
\left\{
\begin{array}{ll}
\nabla \cdot\Gs \nabla u=0 & \mbox{in} \ \ D, \\
\GD u =0 & \mbox{in} \ \ \RR^d \setminus \overline{D}, \\
\ds \Gs \frac{\p u}{\p \nu} \Big|_- = \frac{\p u}{\p \nu} \Big|_+ & \mbox{on} \ \ \p D,\\
u|_- - u|_+=f & \mbox{on} \ \ \partial D,\\
u(x) = O(|x|^{1-d})  & \mbox{as } |x| \rightarrow \infty
\end{array}
\right.
\end{equation}
for any $f\in H_0^{1/2}(\p D)$.
Thus
$\GL_\Gs-\GL^e$ is invertible and $\|\GL_\Gs-\GL^e\|^{-1}\leq C^{-1}$.

Next, we consider the norm of the operator $(\GL_\Gs-\GL^e)^{-1}(\GL_1-\GL^e)$. To get the bounds of this norm,
we only need to consider the bounds of the function
$$
f(x)= \frac{1+x}{a+x} , \quad a>0 \quad x\in [0, \infty).
$$
It can be easily seen that
$$
\min(1,1/a) \leq f(x)\leq \max(1,1/a).
$$
Thus we also get the upper and lower bounds for $(\GL_\Gs-\GL^e)^{-1}(\GL_1-\GL^e)$ and
we come to the conclusion.
\qed

The following stability theorem follows directly from \eqnref{eq:Gsbound} and Theorem 1.1 in \cite{BFR07}.
\begin{theorem}\label{th:stability}
Let $\Gs_1$ and $\Gs_2$ be two different conductivity distribution. Suppose $\Gs_1$ and $\Gs_2$ are H$\ddot{o}$lder continuous, i.e.,
$\Gs_i\in C^s(D)$ with $s>0$. Then
\beq
\|\Gs_1-\Gs_2\|\leq C (\log\|M_{\Gs_1}-M_{\Gs_2}\|)^{-t},
\eeq
where positive constants $C$ and $t$ depend on $\Gs_i$, $i=1,2$ and $s$.
\end{theorem}
We have shown that the stability of the reconstruction of inhomogeneous conductivity by using GPTs can be
actually connected to the stability of that by using the operator $M_{\Gs}$, and it turns out that the
stability result is the same to that by using NtD map. We remark that the EIT (Electrical Impedance Tomography),
to reconstruct $\Gs$ from boundary measurements, is known to be severely ill-posed.
A log-type stability was obtained by Alessandrini \cite{Ales88} and it is
optimal \cite{Mand01}.  A Lipschitz type stability estimate for the values of the conductivity
from the DtN map was proven in \cite{Nach88,SyUh88}.

\section{Reconstruction of GPTs}
From the last section, we see that reconstruction of H$\ddot{o}$lder continuous inhomogeneous conductivity may
have logarithm stability. Uniqueness of determination of $L^\infty$ inhomogeneous conductivity also suggest that
the GPTs can be used for the reconstruction of inhomogeneous conductivity \cite{ADKL12}. However, in
real world applications, usually we can not do the boundary measurements directly hence GPTs or contracted GPTs can not
be calculated directly from \eqnref{eq:conGPTs}. In stead, the common way is to deploy some transmitters and receivers around the target
we are trying to reconstruct. Usually the transmitters and receivers are deployed far away from the target. We suppose that
$D=\Bz+\delta B=\{\Bx=\Bz+\delta\By|\By\in B\}$, where $B$ is a $\mathcal{C}^2$ and bounded domain which has length scale of order one and
$\Bz$ is the center of $D$.
In this section, we shall use the multistatic response (MSR) Matrix to reconstruct the GPTs.
We follow the steps in \cite{ABGJKW} in getting the GPTs. Suppose $\{\Bx_t\}_{t=1}^{N_r}$ and $\{\Bx_s\}_{s=1}^{N_s}$ are a set of
electric potential point detectors and electric point sources. We suppose that the number and location of the point detectors
coincide with the number and location of the point sources, i.e., $N_t=N_s=N$. Then the MSR matrix, denoted by $\BV$ is an $N\times N$
matrix whose $ts$ entry $V_{ts}$ is defined by
\beq\label{msr}
V_{ts}=u_s(\Bx_t)-\Gamma_s(\Bx_t),\quad t,s=1,2,\ldots, N.
\eeq
Here $\Gamma_s(\Bx)=\Gamma(\Bx-\Bx_s)$. $u_s(\Bx)$ is the solution to the following transmission problem
\beq\label{eq:msrtran}
\left\{
\begin{array}{ll}
\nabla\cdot \Gs \nabla u_s(\Bx)= \delta_{\Bx_s}, & \Bx\in \RR^2\setminus \p D\\
u_s(\Bx)|_+=u_s(\Bx)|_- & \Bx\in \p D\\
\frac{\p u_s}{\p \bm\nu}(\Bx)\Big|_+= \Gs \frac{\p u_s}{\p \bm\nu}(\Bx)\Big|_-  & \Bx\in \p D\\
u_s(\Bx)-\Gamma_s(\Bx)=O(|\Bx|^{-1}) & \mbox{as}\quad |\Bx-\Bx_s|\rightarrow \infty.
\end{array}
\right.
\eeq
The solution to \eqnref{eq:msrtran} is unique and satisfies
\beq
u_s(\Bx)- \Gamma_s(\Bx)= \Scal_D[\phi_s](\Bx)=\int_{\p D} \Gamma(\Bx-\By) M_{\Gs}\Big[\frac{\p \Gamma_s}{\p \bm\nu}\Big|_{\p D}\Big]ds(\By)
, \quad \Bx\in \RR^2\setminus \overline{D}.
\eeq
For any $\By\in \p D$ and $\Bz$ away from $\Bx$ we have the $K$-th order Taylor expansion for $\Gamma(\Bx-\By)$
$$
\Gamma(\Bx-\By)=\Gamma(\Bx-\Bz-(\By-\Bz))=\sum_{|\alpha|=1}^{K} \frac{(-1)^{|\alpha|}}{\alpha!}\p^{\alpha} \Gamma(\Bx-\Bz) (\By-\Bz)^{\alpha}+ e_K,
$$
where $e_K$ is the truncated error. The MSR matrix thus has the following approximation
\beq\label{eq:MSR2GPT}
V_{ts}=\sum_{|\alpha|=1,|\beta|=1}^{K} \frac{(-1)^{|\beta|}}{\alpha!\beta!}\p^{\beta} \Gamma(\Bx_t-\Bz)Q_{\alpha \beta}\p^{\alpha} \Gamma(\Bz-\Bx_s) +{E}_{ts}
\eeq
where
$$
Q_{\alpha \beta}=\int_{\p D} (\By-\Bz)^{\beta} M_{\Gs}\Big[\frac{\p (\cdot-\Bz)^{\alpha}}{\p \bm\nu}\Big|_{\p D}\Big](\By)ds(\By).
$$
Note that the zeroth order expansion varnish because $M_{\Gs}$ maps from $H_0^{-1/2}$ to $H_0^{-1/2}$ (which makes the term $|\beta|=0$ varnish) and there is a normal derivative $\p/{\p \bm\nu}$ (which makes the term $|\alpha|=0$ varnish). Since $\GL^e=I_0\Scal_D(1/2I+\Kcal_D^*)^{-1}$ is invariant under the translation of the domain $D$, by using the fact that $\By-\Bz\in \p (\delta B)$ if $\By\in \p D$ we have
$$
Q_{\alpha \beta}=\int_{\p (\delta B)} \tilde{\By}^{\beta} M_{\tilde{\Gs}}\Big[\frac{\p \tilde{\Bx}^{\alpha}}{\p \bm\nu}\Big|_{\p (\delta B)}\Big]ds(\tilde{\By})
= M_{\alpha\beta}(\tilde{\Gs})
$$
where $\tilde{\Gs}(\tilde{\Bx}):=\Gs(\Bx-\Bz)$ for $\tilde{\Bx}\in \delta B$. By using the optimization method we can thus
reconstruct the GPTs of the shifted domain $\delta B$. Since it is more convenient to use the contracted GPTs to do the
reconstruction, we shall present the optimization method to get the CGPTs. To do this, we should first expand the fundamental
solution in a different way, or using harmonic expansion (see, e.g., \cite{ADKL12,nedelec})
 \begin{align}\label{Gammaexp1}
 \GG(\Bx-\By)
 =-\sum_{n=1}^{\infty}\frac{1}{2\pi n}\left[\frac{\cos n\theta_x}{r_x^n}r_y^n\cos n\theta_y
 +\frac{\sin n\theta_x}{r_x^n}r_y^n\sin n\theta_y\right]+ C,  \quad d=2
 \end{align}
\begin{align}
\Gamma(\Bx-\By)= - \sum_{\ell=0}^\infty \sum_{k=-\ell}^{\ell} \frac{1}{2\ell+1}Y_\ell^k(\theta_x, \varphi_x) \overline{Y_\ell^k(\theta_y, \varphi_y)} \,\frac{r_y^l}{r_x^{l+1}}, \quad d=3 \label{Gammaexp2}
\end{align}
By using the harmonic expansions (d=2) we then have
\beq\label{eq:MSR2GPT}
V_{ts}=\sum_{m,n=1}^{K} \bm{A}_{tm}^T\bm{M}_{mn}(\tilde{\Gs})\bm{A}_{ns} +{E}_{ts}
\eeq
where
$$\bm{A}_{tm}:=\frac{1}{2\pi m r_{t}^m}(\cos m\theta_{t}, \sin m\theta_{t})^T$$
and $T$ means the transpose of a vector and $\bm{M}_{mn}(\tilde{\Gs})$ is a two-by-two matrix which has the form
\begin{align*}
\bm{M}_{mn}(\tilde{\Gs}):=\left[
\begin{array}{ll}
M_{mn}^{cc}(\tilde{\Gs}) & M_{mn}^{cs}(\tilde{\Gs})\\
M_{mn}^{sc}(\tilde{\Gs}) & M_{mn}^{ss}(\tilde{\Gs})
\end{array}
\right].
\end{align*}
If $d=3$ then it has different dimensions for different $m$ and $n$. The result can be similarly got by using
the expansion \eqnref{Gammaexp2} (see \cite{ACKW12}).
Denote by $\bm{M}(\tilde{\Gs})=(\bm{M}_{mn}(\tilde{\Gs}))_{m,n\leq K}$ the matrix combined by the GPTs of order lower than $K$. Without making ambiguity sometimes we
omit the dependence on $\Gs$ and write $\bm{M}$ instead. The dimension of $\bm{M}$ depends
on the dimension of the space ($d$) and the number $K$. If $d=2$ then $\bm{M}$ is a $2K$-by-$2K$ matrix and if $d=3$ then $\bm{M}$
is a $K(K+2)$-by-$K(K+2)$ matrix. We can consider $\bm{M}$ as a lower frequency approximation of the operator $M_{\Gs}$.
Define $\bm{A}:=(\bm{A}_{mt})_{m\leq K, t \leq N}$, a $2K$-by-$N$ matrix and $\bm{E}:=(\bm{E}_{ts})_{t, s\leq N}$, a $N$-by-$N$ matrix then we have
$$
\bm{V}=\bm{A}^T \bm{M} \bm{A}+\bm{E}
$$
The reconstruction of the CGPTs is achieved then by using
the following optimization method
$$
\BM^{est}:=\min_{\BM\in\RR^{2N\times 2N}}\|\BV-\bm{A}^T\BM\bm{A}\|^2_F
$$
where $\|\cdot\|_F$ is a matrix norm which can be Frobenius, etc. In general, in order to reconstruct CGPTs of order
$K$, $N$ should be taken large enough such that $2K<N$ (for $d=2$). Observe from \eqnref{Gammaexp1} and
expression of MSR matrix that the contribution of a CGPT decays as its order grows. Consequently, the inverse procedure
may not be stable for higher order CGPTs. Stability analysis is thus required.
When the MSR matrix is measured with
noise, stability results can also be obtained follow similar steps in \cite{ABGJKW}, where excellent analysis and comments
are given.

\section{Optimization method for reconstruction of the inhomogeneous conductivities}
It follows from \eqnref{eq:Gsbound} that we can use the information of the operator $M_{\Gs}$ in stead of the NtD map $\GL_{\Gs}$ to reconstruct the conductivity $\Gs$.
Numerically, an effective way is to approach the eigenvalues of $M_{\Gs}$ with respect to low frequency eigenfunctions. Since $M_{mn}(\Gs)$
plays the role of the eigenvalues of $M_{\Gs}$ when $h_m$ and $h_n$ are chosen appropriately, reconstruction of $M_{mn}$ for some lower
numbers $m,n \in \NN$ would be a natural way to reconstruct the conductivity $\Gs$. Higher order GPTs are known to be quite unstable and,
in the mean time, they require much more computations as the order of GPTs increase.

\subsection{Finite approximation of the operator}
We use a least-square approach (see, for
instance, \cite{hanke}) for the reconstruction of $\Gs$. Let $\Gs^*$ be the exact (target)
conductivity (in two dimensions) and let $y_{mn}:=M_{mn}(\Gs^*)$. One optimization approach (cf. \cite{ADKL12}) is to
minimize the following discrepancy functional
\beq
\label{eq:lsq1} S_1(\Gs)=\frac{1}{2}\sum_{m,n \le N}
\omega_{mn}|y_{mn}-M_{mn}(\Gs)|^2
\eeq
for some well-chosen weights $\omega_{mn}$. The weights $\omega_{mn}$ plays an important role in convergence rate of
minimizing the functional. In general, $\omega_{mn}$ should depend on the eigenvalues of the matrix $\bm{M}$. It is then natural
to introduce another functional (see \cite{siims})
\beq
\label{eq:lsq2}
S_2(\Gs)=\frac{1}{2}\sum_{l=1}^N \sum_{l'=1}^N \omega_l(\Gl_0^{(l)})\omega_{l'}(\Gl^{(l')})|\la(\bm{Y}-\bm{M})\bm{v}_0^{(l)},\bm{v}^{(l')}\ra|^2
\eeq
where $\bm{Y}:=(y_{mn})$. $\Gl_0^{(l')}$ and $\Gl^{(l)}$ are the eigenvalues of $\bm{Y}$ and $\bm{M}$, respectively. $\bm{v}_0^{(l)}$ and $\bm{v}^{(l')}$
are the corresponding eigenvectors of $\bm{Y}$ and $\bm{M}$, respectively. It is clearly seen that if $\bm{Y}$ and $\bm{M}$ are diagonal matrix then
\eqnref{eq:lsq2} and \eqnref{eq:lsq1} are the same formula and if $\bm{Y}$ is a small perturbation of a diagonal matrix, then one
approximate another.

\subsection{Calculation of Fr\'echet derivatives and reconstruction algorithms}

The Fr\'echet derivative in the direction of
$\Gg$, a function in $D$, is define to be
 $$
 M'_{mn}(\Gs)[\Gg] := \lim_{\ep \to 0} \frac{M_{mn}(\sigma+ \ep \Gg)-M_{mn}(\sigma)}{\ep}.
 $$
From \cite{ADKL12} we have
\beq\label{eq:Fre_op}
M'_{mn}(\Gs)[\gamma]= \int_{D} \gamma \na u_n \cdot \na u_m d\Bx,
\eeq where $u_n$ and $u_m$ are the solutions of
\begin{equation} \label{defh}
\left\{
\begin{array}{ll}
\na \cdot (\Gs \na u)=0 & \mbox{in }  D, \\
\nm \ds \Gs \frac{\p u}{\p
\bm\nu}=(\GL_{\Gs}-\GL^e)^{-1}(\GL_{1}-\GL^e)[\na h \cdot \bm\nu] &
\mbox{on } \p D,
\end{array}
\right.
\end{equation}
with $h= h_n$ and $h=h_m$, respectively.
Note that if $M_{mn}(\sigma)$ is one of the other contracted GPTs,
then $h$ should be changed accordingly.
One can easily see that the adjoint $M'_{mn}(\Gs)^*$ of
$M'_{mn}(\Gs)$ is given by
 \beq\label{eq:Fre_op1}
 M'_{mn}(\Gs)^*[c]=c \na u_m \cdot \na u_n, \quad c \in \RR
 \eeq
and we sometimes write $M'_{mn}(\Gs)^*=\na u_m \cdot \na u_n$.

Denote by $(\bm{M}'(\Gs))_{mn}:=M_{mn}'(\Gs)$ the Fr\'echet derivative of the matrix $\bm{M}$ then we have
the Fr\'echet derivative of the functional $S_2(\Gs)$
$$
S_2'(\Gs)[\gamma]=-\sum_{l=1}^N \sum_{l'=1}^N \omega_l(\Gl_0^{(l)})\omega_{l'}(\Gl^{(l')})\la(\bm{Y}-\bm{M})\bm{v}_0^{(l)},\bm{v}^{(l')}\ra
\la\bm{M}'(\Gs)[\gamma]\bm{v}_0^{(l)},\bm{v}^{(l')}\ra
$$
and so the adjoint
$$
S_2'(\Gs)^*[c]=-\sum_{l=1}^N \sum_{l'=1}^N \omega_l(\Gl_0^{(l)})\omega_{l'}(\Gl^{(l')})\la(\bm{Y}-\bm{M})\bm{v}_0^{(l)},\bm{v}^{(l')}\ra
\la\bm{M}'(\Gs)^*[c]\bm{v}_0^{(l)},\bm{v}^{(l')}\ra.
$$
Next, we consider the algorithms for minimizing the functional $S_2(\Gs)$. Basically, there are two classical iteration methods to minimize
\eqnref{eq:lsq2}: gradient descent method (Landweber) and Newton method. We mention that both methods have their own merit, Landweber method
is relatively more stable and always a descent method but very slow, while Newton method is quite fast but requires a very good initial guess.

The gradient descent procedure to solve the least-square problem \eqnref{eq:lsq1} and \eqnref{eq:lsq2} reads
\beq \label{eq:Land1}
\Gs_{k+1}= \Gs_k - tS_j'(\Gs_k)^*, \quad j=1,2\eeq
where $t$ is a small positive parameter to ensure that $S_j(\Gs_{k+1})<S_j(\Gs_k)$. For step chosen, we can apply the classical
Armijo's rule. On the other hand the Newton method reads
\beq \label{eq:Newt}
\Gs_{k+1}=\Gs_k - (S_j'(\Gs_k)^*S_j'(\Gs_k))^{\dag}S_j'(\Gs_k)^*[S_j(\Gs_k)]  \quad j=1,2
\eeq
where $\dag$ means the pseudo-inverse. The newton method requires to get the pseudo-inverse of $S_2'(\Gs)$ which makes the procedure
unstable when the problem is ill-posed. We shall follow the same steps in \cite{siims} to get the pseudo-inverse. Let $\gamma$ belong to
the vector spaces spanned by $\{\psi_p\}$
$$
\{\psi_p\}=\{\la\bm{M}'(\Gs)^*\bm{v}_0^{(l)},\bm{v}^{(l')}\ra\}
$$
then by \eqnref{eq:Fre_op} and \eqnref{eq:Fre_op1} we have
$$
S_2'(\Gs)^*[c]=c\sum_p S_2'(\Gs)[\psi_p]\psi_p
$$
and by using the identities of pseudo-inverse
$$
(S_2'(\Gs))^{\dag}=(S_2'(\Gs)^*S_2'(\Gs))^{\dag}S_2'(\Gs)^*=S_2'(\Gs)^*(S_2'(\Gs)S_2'(\Gs)^*)^{\dag}
$$
we obtain
\beq \label{eq:pseu_in}
(S_2'(\Gs)^*S_2'(\Gs))^{\dag}S_2'(\Gs)^*[S_2(\Gs)]= \frac{S_2(\Gs)}{\sum_p |S_2'(\Gs)[\psi_p]|^2} \sum_p S_2'(\Gs)[\psi_p]\psi_p.
\eeq
Note that the larger $N$ is, the more components of $\Gg$ shall be reconstructed and so better resolution may be
obtained.
We shall consider the resolution of the reconstruction. The notion of the resolution in solving inverse conductivity
 problem was first introduced in \cite{procams}. The definition of resolution is inherited form
 classical Rayleigh resolution formula for active array imaging \cite{BPTs03,BoWo99}.
 For simplicity we assume that the exact conductivity distribution $\Gs^*$ in
a disk $B$ is a perturbed constant conductivity $c$ in $B$. We denote $\Gs_0:=c$.

\subsubsection{First functional}
In this subsection, we use \eqnref{eq:lsq1} to reconstruct the conductivity distribution.
On the one hand, by using \eqnref{eq:Land1} we get the first perturbation
$$\gamma:=\Gs_1-\Gs_0=-tS_1'(\Gs_0)^*=t\sum_{m,n\leq N}\omega_{mn}\nabla u_m\cdot \nabla u_n (y_{mn}-M_{mn}(\Gs_0)).$$
Note that $\Gs_0$ is constant in a disk $B$, we compute directly $u_m$ and $M_{mn}(\Gs_0)$
$$
u_m=2/(c+1) h_m \quad \mbox{and} \quad
M_{m,n}(\Gs_0)=\left\{
\begin{array}{ll}
0 & m\neq n \\
2\pi m (c-1)/(c+1)r_0^{2m} & m=n
\end{array}
\right.
$$
Thus
\beq\label{eq:reso1}
\gamma(r,\theta)=8/(c+1)^2t\sum_{m,n\leq N}mn\omega_{mn} r^{m+n-2} e^{i(m-n)\theta}(y_{mn}-M_{mn}(\Gs_0))
\eeq
where we choose $h_m=e^{im\theta}$. If furthermore the exact conductivity $\Gs^*$ is radial symmetric and (analytic) can be written
$$
\Gs^*(r)= \sum_{m=0}^{\infty} a_n r^n
$$  then $y_{mn}=0$ for $m\neq n$ and we see that any polynomial of order less than $2N-2$ can be reconstructed by one step if the
weights $\omega_{mn}$ is chosen optimally. However, polynomials of order larger than $2N-2$ can not be reconstructed by using GPTs
of order up to $N$.

On the other hand, By using \eqnref{eq:Newt} we get the first perturbation
\beq\label{eq:reso2}
\gamma(r,\theta)=\sum_{m,n\leq N} C_{mn} r^{m+n-2}e^{i(m-n)\theta}
\eeq
where the coefficients $C_{mn}$ are
\begin{align*}
C_{mn}=\frac{mn\omega_{mn}\epsilon_{mn}\sum_{m'n'\leq N}\omega_{m'n'}\epsilon_{m'n'}^2}{4\pi\sum_{m,n,m',n'\leq N}\omega_{mnm'n'}
\epsilon_{mn}\epsilon_{m'n'}r_0^{m+n+m'+n'-2}\delta_{m+m'n+n'}}
\end{align*}
where $\delta_{mn}$ is the Kronecker delta function and
$$
\omega_{mnm'n'}=\frac{mnm'n'}{m+n+m'+n'-2}\omega_{mn}\omega_{m'n'}\quad \mbox{and} \quad \epsilon_{mn}=y_{mn}-M_{mn}(\Gs_0).
$$
We can thus get the same result by choosing appropriately the weights $\omega_{mn}$.
For general conductivity, by \eqnref{eq:reso1} and \eqnref{eq:reso2}, we have the following conclusions:
\begin{enumerate}
  \item[i] The larger $N$ is the better angular resolution in reconstructing of the conductivity can be got. In the meantime,
  suppose the conductivity is smooth enough in $B$ and has the Taylor expansion
  $$
  \Gs^*=\sum_{|\alpha|=0}^\infty \frac{1}{\alpha!}\p^{\alpha}\Gs^*(0) r^{|\alpha|}\cos^{\alpha_1}\theta \sin^{\alpha_2}\theta.
  $$
  We thus conclude that only the coefficients of expansions with order less than $2N-2$ can be reconstructed by using orders of CGPTs which
  are less than $N$.
  \item[ii] Lower order GPTs only contains the low-frequency information on reconstructing $\gamma$, while
  higher order GPTs contains both low-frequency and high-frequency information. And so reconstruction of $\gamma$
  by higher order GPTs is quite unstable.
  \item[iii] Reconstruction of $\gamma$ near the origin
(especially $r=0$) is more sensitive to noise than near the boundary of $B$.
This is in accordance with \cite{procams, gunther}.
  \item[iv] From the stability result in section 3, the convergence of the H$\ddot{o}$lder continuous conductivity is as slow as
  the inverse of logarithm of the convergence of the GPTs. Thus by \eqnref{eq:reso1} and \eqnref{eq:reso2}, we conclude that the
  convergence rate is extremely slow for approximating the higher order GPTs.
\end{enumerate}
\subsubsection{Second functional}
For the functional $S_2(\Gs)$ if we use \eqnref{eq:Land1} then we get the lower order of the first perturbation
$$
\gamma\sim8/(c+1)^2t\sum_{l,l'\leq N} \omega_{ll'} ll' r^{l+l'-2} e^{i(l-l')\theta}(y_{ll'}-M_{ll'}(\Gs_0))
$$
where $\omega_{ll'}=\omega_l(\Gl_0^{(l)})\omega_{l'}(\Gl^{(l')})$. It is seen that if the exact conductivity is a small
perturbation of a constant conductivity in a disk then both functionals converge to each other. We can also get the
estimate of the perturbation if we use the Newton method by using similar methods.

\subsection{Regularization due to measurement noise}
In this section, we consider the reconstruction of the conductivity due to measurement noise.
In real applications, the GPTs are usually reconstructed from the MSR matrix. As is done before, the reconstructed GPTs
are not accurately given. For the sake of convenience, we denote by $y_{mn}^{\delta}$ the reconstructed GPTs and $\bm{Y}^{\delta}$
the matrix which satisfies
$$
\|\bm{Y}^{\delta}-\bm{Y}\|_F \leq \delta
$$
In this case, the least square functional
\eqnref{eq:lsq1} and \eqnref{eq:lsq2} can not be used directly to reconstruct the inhomogeneous conductivity due to the
severely ill-posedness of this problem. Regularization method is then required. For analysis of regularization methods,
refer to \cite{hanke,Loui96,Loui99}.  According to \eqnref{eq:lsq1},
we introduce the regularized least square functional as follows
\beq
\label{eq:lsq3} S_3(\Gs)=\frac{1}{2}\sum_{m=1}^N \sum_{n=1}^N \omega_{mn}|y_{mn}^{\delta}-M_{mn}(\Gs)|^2 + q\|\Gs-\Gs_0\|_p^2
\eeq
and according to \eqnref{eq:lsq2} we introduce
\beq
\label{eq:lsq4} S_4(\Gs)=\frac{1}{2}\sum_{l=1}^N \sum_{l'=1}^N \omega_l(\Gl_0^{(l)})\omega_{l'}(\Gl^{(l')})
|\la(\bm{Y}^{\delta}-\bm{M})\bm{v}_0^{(l)},\bm{v}^{(l')}\ra|^2 + q\|\Gs-\Gs_0\|_p^2
\eeq
where $\Gs_0$ is an \emph{a priori} information on the exact conductivity $\Gs^*$ and $q$ is a small positive parameter
called regularization factor. $\|\cdot\|_p$ is a general norm which can be $L^2$ or $L^{\infty}$ and so on.

To simplify the analysis, we only consider the functional \eqnref{eq:lsq3}. It can be proved (see e.g. \cite{hanke})
that if the regularization factor $q$ which depends on the noise level $\delta$ satisfies some basic convergence rules
($q\rightarrow 0$ and $\delta^2/q\rightarrow 0$ as $\delta\rightarrow 0$) then the minimizer to \eqnref{eq:lsq3}
converges to the $\Gs_0$-minimum-norm solution. If $L^2$ norm is chosen in the penalty term $\|\Gs-\Gs_0\|_p^2$, then
we can directly use the following Landweber iteration to get the minimizer of \eqnref{eq:lsq3},
$$
\Gs_{k+1}=\Gs_k+t\sum_{m,n\leq N}\omega_{mn}M_{mn}'(\Gs_k)^*[y_{mn}^{\delta}-M_{mn}(\Gs_k)]-2tq(\Gs_k-\Gs_0).
$$
Since Landweber iteration is itself a regularization method, in the numerical experiments we let $q=0$.
We mention that the iteration index plays the role of the regularization parameter, and the stopping criterion is the
counterpart of the parameter choice rule in continuous regularization methods \cite{hanke}. Consequently, appropriate stopping rules
are essential in getting a good approximation result. We use the widely-used Morozov discrepancy principle as the
\emph{a posteriori} stopping rule \cite{MO:1966}. The initial estimation $\Gs_0$ plays an important role in
reconstruction of the conductivity and should be chosen such that it is "near" the exact conductivity.
Since in numerical implementation, not all information of GPTs can be used for reconstruction, in fact only
several lower order GPTs are used as we shall see, the uniqueness of the inverse problem in reconstruction of
conductivity is not ensured. Thus different initial approximation may produce different reconstruction result
($\Gs_0$-minimum-norm solution), especially for discontinuous conductivity, where iteration may go to
the reconstruction of smooth conductivity which has the same lower order GPTs.

\section{Numerical experiments}
In this section, we consider the reconstruction of four main different types of conductivities, that is, antisymmetric
conductivity, general $C^{\infty}$ conductivity and H$\ddot{\mbox{o}}$lder continuous conductivity and
discontinuous conductivity. The four different types of
conductivities are as follows:
\begin{enumerate}
  \item[1.] $\Gs=x^3+y^3+4.0$
  \item[2.] $\Gs=bx^3+ay^5+y^2+2.0$
  \item[3.] $\Gs=x^3+y^5+(y-0.5)^{0.4}+3.0$
  \item[4.] $\Gs=x^3+y^5+\chi(x^2+y^2<0.25)+2.0$
\end{enumerate}
where $a$ and $b$ are given constants. The inclusion $D$ is supposed to be a unit disk.
In the numerical implementation, the GPTs are computed up to order five or six from the boundary measurements. The boundary
measurements are calculated by solving the corresponding partial differential equations which may produce some measurement
errors. We can also reconstruct the GPTs from MSR matrix, but in this paper we only focus on the reconstruction of the
conductivities by using GPTs thus we skip that part.
We use a very fine mesh to compute $y_{mn}, 1\leq m, n \leq 6$. To stably and accurately reconstruct the conductivity
distribution, we use a recursive approach proposed in \cite{AKLZ12} (see also \cite{siims, AGKLY11, bao, borcea2}).
We first minimize the discrepancy between the first contracted GPTs for $1\leq m, n \leq l$, where $l< 6$ by using a coarse mesh.
Then we use the result as an initial guess for the minimization between the GPTs for $1\leq m, n \leq l+1$.  This corresponds
to choosing appropriately the weights $\omega_{mn}$ in \eqnref{eq:lsq1}. In addition, the weights $\omega_{mn}$ also play the
role for the convergence of the iteration method \eqnref{eq:Land1} (also called Landweber iteration).  Moreover, we refine the mesh
used to compute the reconstructed conductivity distribution every time we increase the number of used contracted GPTs in the
discrepancy functional.

Let $k$ be the iteration step and let $\varepsilon_{M}$ and $\varepsilon_{\Gs}$ be discrepancies of GPTs and the conductivities, {\it i.e.},
 \beq
 \varepsilon_{M}:= \sum_{m,n \le N}(y_{mn}-M_{mn}(\Gs_{k}))^2, \quad y_{mn}= M_{mn}(\Gs^*),
 \eeq
($N$ represents the order of GPTs used) and
 \beq
 \varepsilon_{\Gs}:=\frac{\int_{ D}
(\Gs_{k}-\Gs^*)^2}{\int_{ D} (\Gs^*)^2}.
 \eeq
Figure \ref{fig-result}--\ref{fig3-result} are the reconstructed first three conductivity distributions ($a=1.0, b=1.0$)
using contracted GPTs with order up to $N=5$, $N=6$ and $N=5$, respectively. It can be seen from these figures that the reconstruction of antisymmetric conductivity
has the best accuracy, or is easiest, among the three different types of conductivities. And we have very accurate approximation for the antisymmetric
conductivity and the general $C^{\infty}$ conductivity. This is because that the regularities of the first two conductivities to be reconstructed are much better than the third conductivity to be reconstructed, which also involve less measurement errors than the third one in computing the partial differential equations. In addition, Figure \ref{fig2-result} shows that there is little improvement in the performance of reconstruction by $N=6$ GPTs than by $N=5$ GPTs, which is due to the reason that higher order GPTs is very unstable to reconstruct and it mainly contains the high frequency parts of the
conductivity. In Figure \ref{fig3-result} we see that the inaccuracy of reconstruction mainly occurs near the discontinuous region of the conductivity
($y=0.5$). Figures \ref{fig-result1}-\ref{fig3-result1} are the reconstruction history of $\varepsilon_{M}$ and $\varepsilon_{\Gs}$. The sudden jump in the figures happen when the number of GPTs used for reconstruction changes from $N$ to $N+1$. From these figures we can find that the convergence rate
decays quite fast after few steps in each reconstruction with fixed number of GPTs.

\begin{figure}[h]
 \begin{center}
   \end{center}
   \caption{Reconstructed conductivity distribution ($\Gs=x^3+y^3+4.0$). The figures from left to right,
   top to bottom are respectively: GPTs order up to one used for approximation,
   GPTs order up to two used for approximation, GPTs order up to three used for approximation,
    GPTs order up to four used for approximation, GPTs order up to five used for approximation and the exact conductivity.\label{fig-result}}
 \end{figure}

\begin{figure}[ht]
\begin{center}
  \end{center}
  \caption{Reconstructed conductivity distribution ($\Gs=x^3+y^5+y^2+2.0$).  The order of GPTs used in turn
  from the upper-left to down-right is from up to order one to up to order six. The bottom figure is the exact conductivity.\label{fig2-result}}
\end{figure}

 \begin{figure}[h]
 \begin{center}
   \end{center}
   \caption{Reconstructed conductivity distribution ($\Gs=x^3+y^5+(y-0.5)^{0.4}+3.0$). The down-right figure is the exact conductivity
   and others are orders of GPTs used for reconstruction from one to five.\label{fig3-result}}
 \end{figure}

%
%
%

\begin{figure}[ht]
\begin{center}
  \includegraphics[width=2.7in]{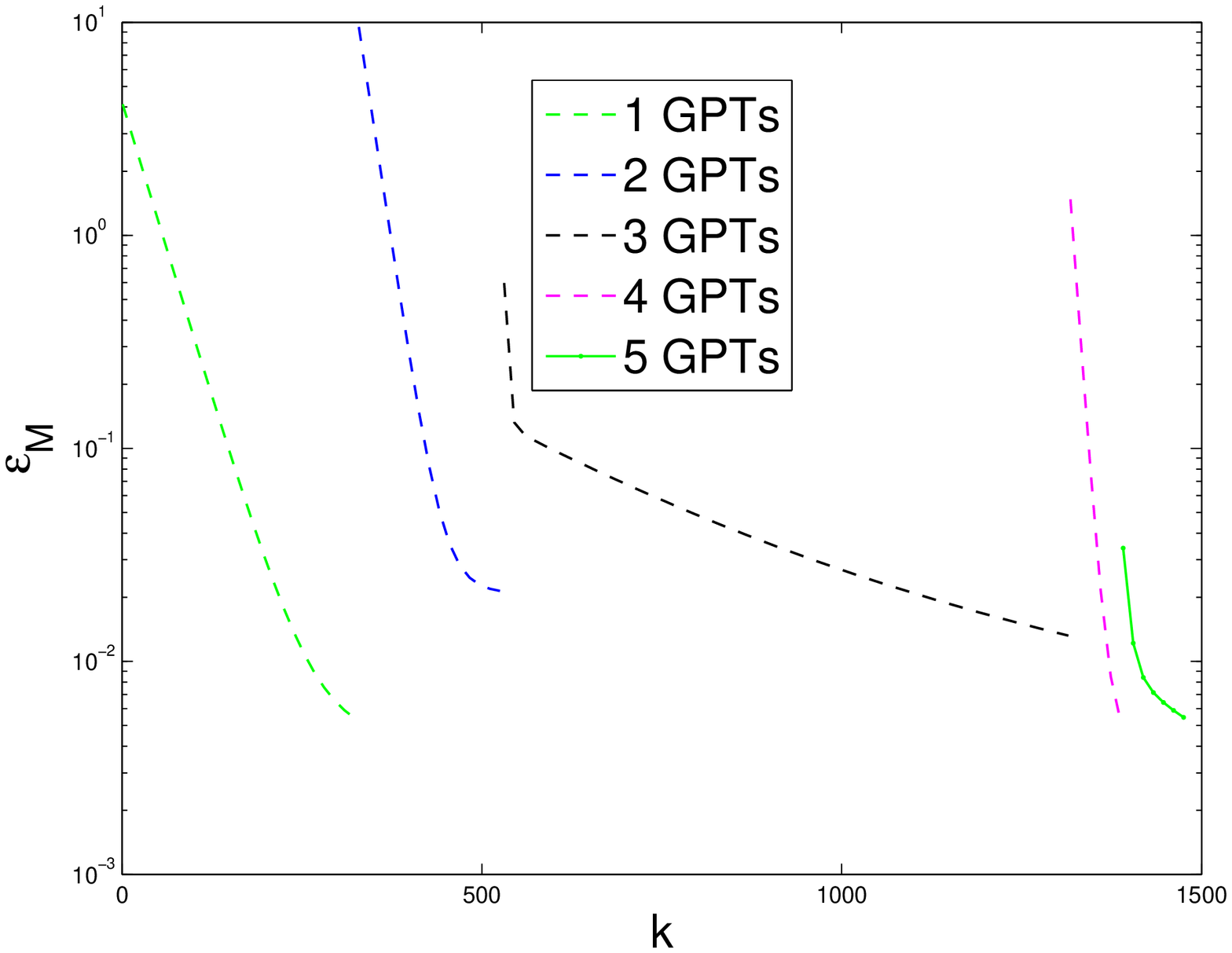}
   \includegraphics[width=2.7in]{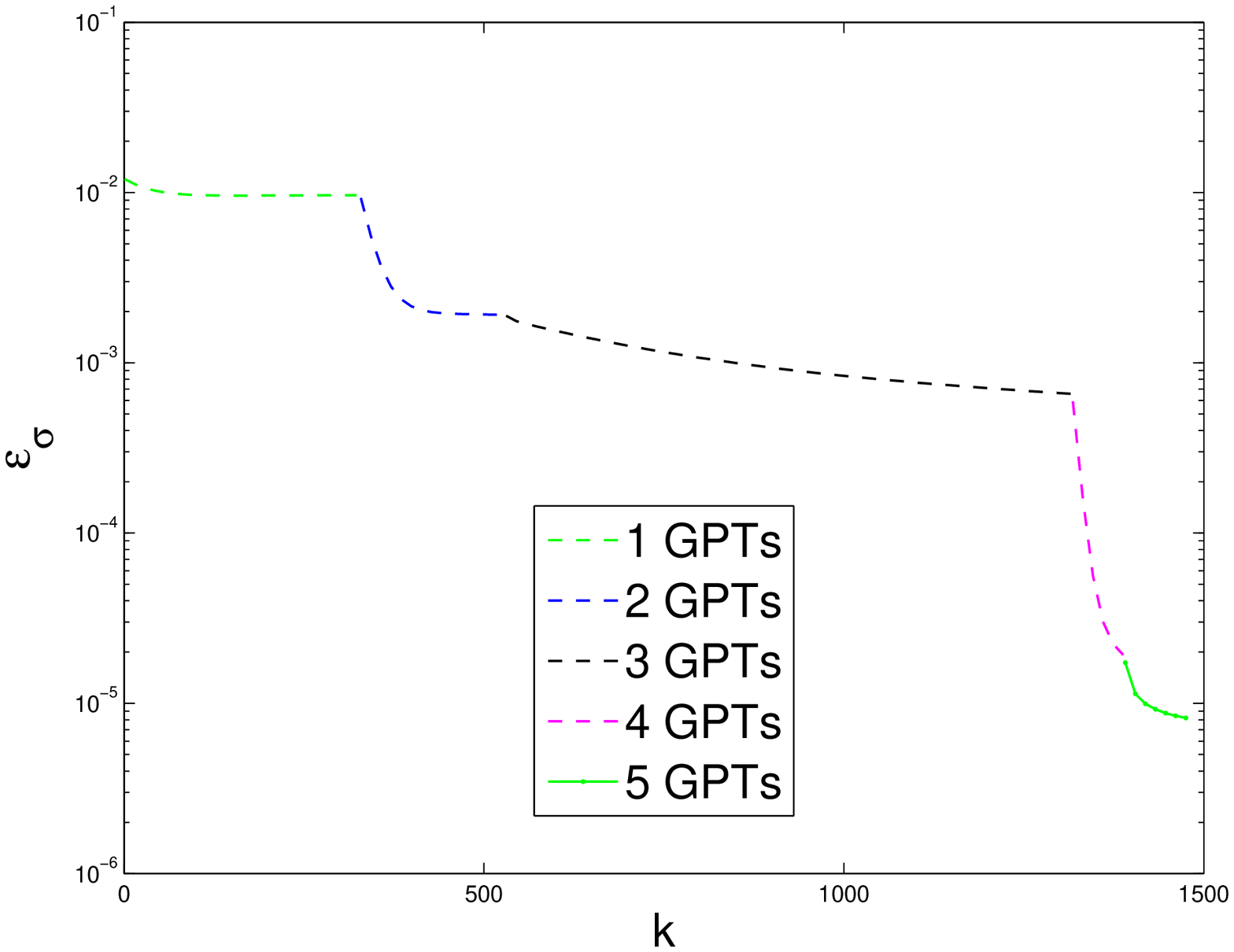}
  \end{center}
  \caption{The convergence history of $\varepsilon_{M}$ and $\varepsilon_{\Gs}$, where ${\rm k}$ is the number of iterations. $\Gs=x^3+y^3+4.0$. \label{fig-result1}}
\end{figure}

\begin{figure}[ht]
\begin{center}
  \includegraphics[width=2.7in]{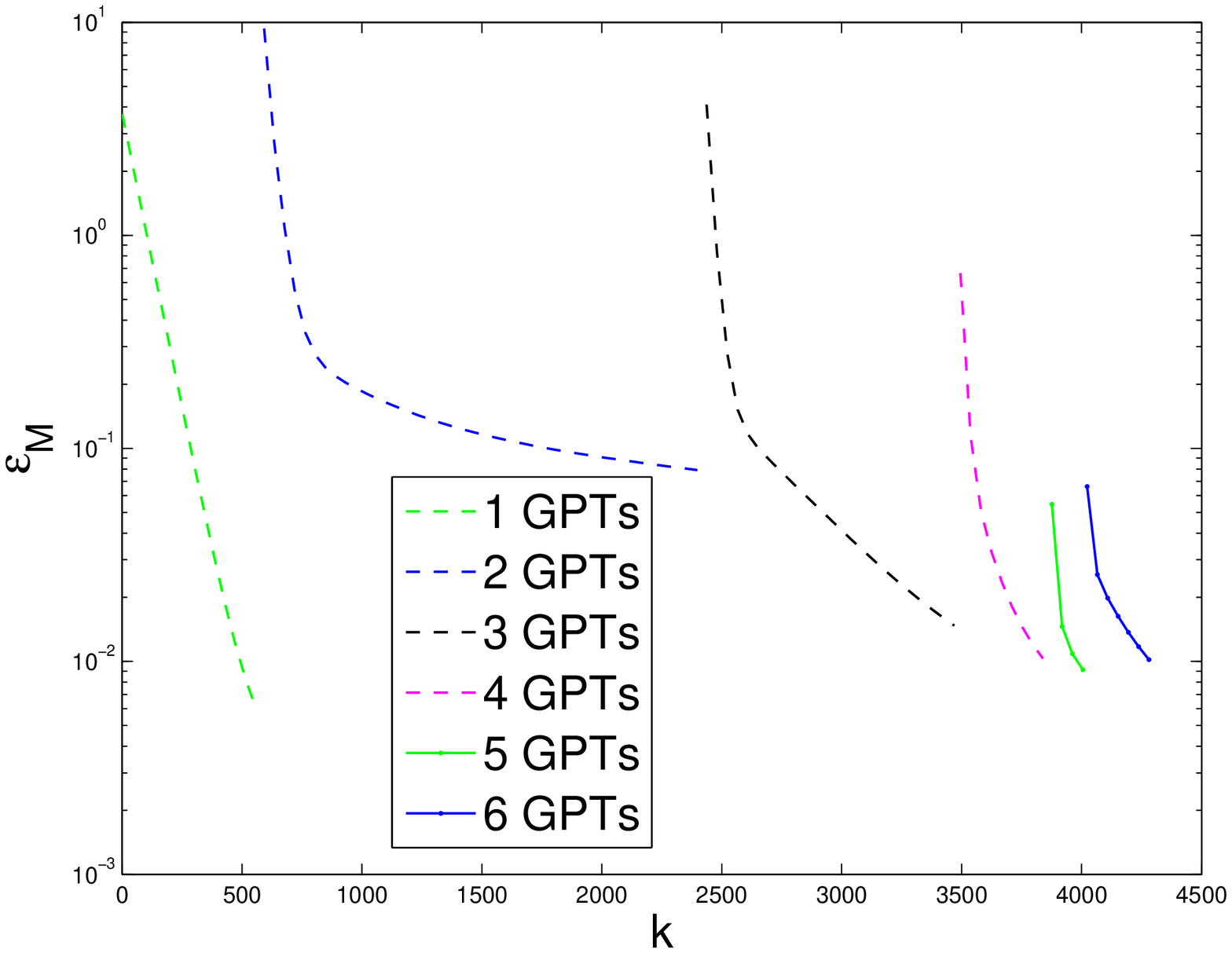}
   \includegraphics[width=2.7in]{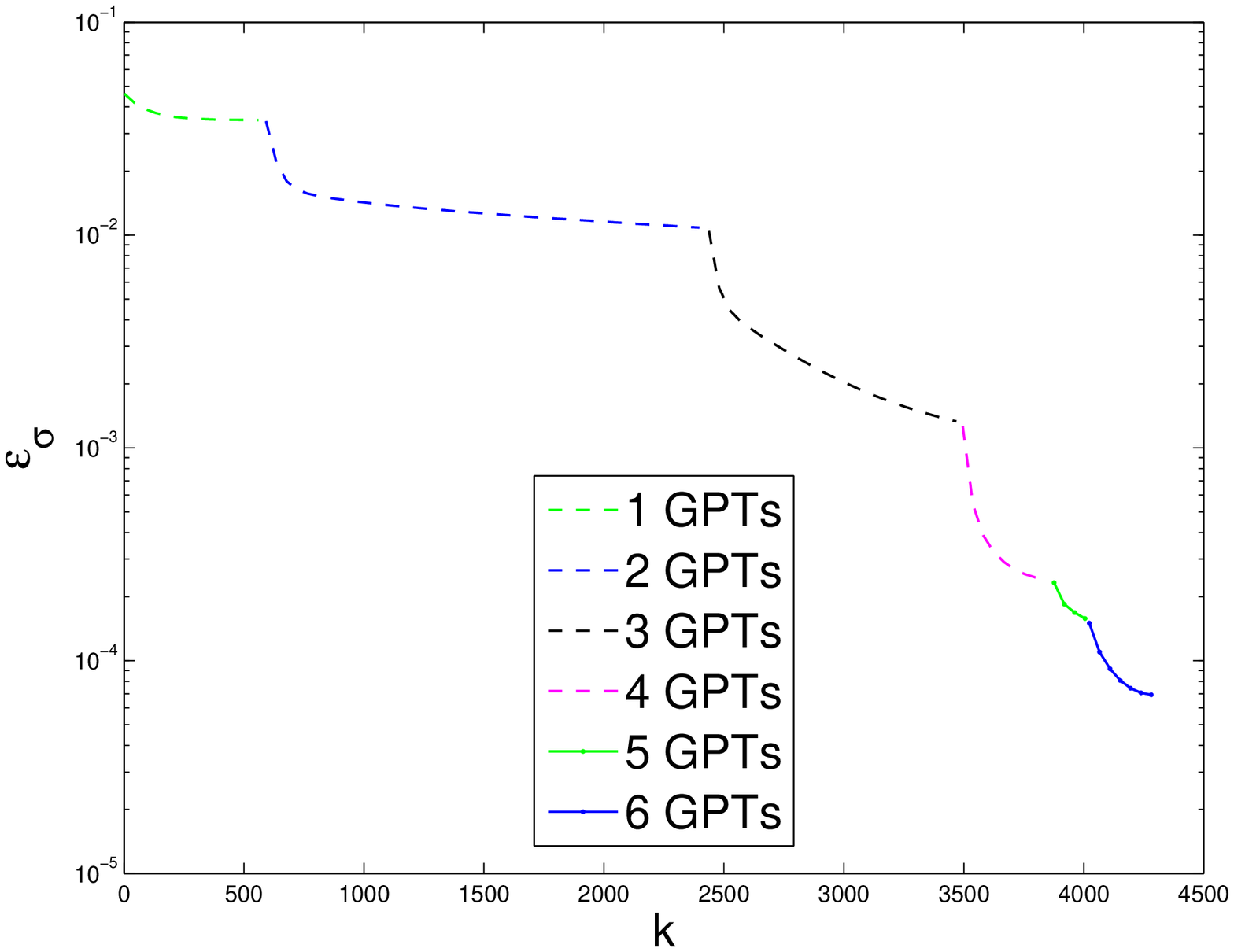}
  \end{center}
  \caption{The convergence history of $\varepsilon_{M}$ and $\varepsilon_{\Gs}$, where ${\rm k}$ is the number of iterations. $\Gs=x^3+y^5+y^2+2.0$. \label{fig2-result1}}
\end{figure}

\begin{figure}[ht]
\begin{center}
  \includegraphics[width=2.7in]{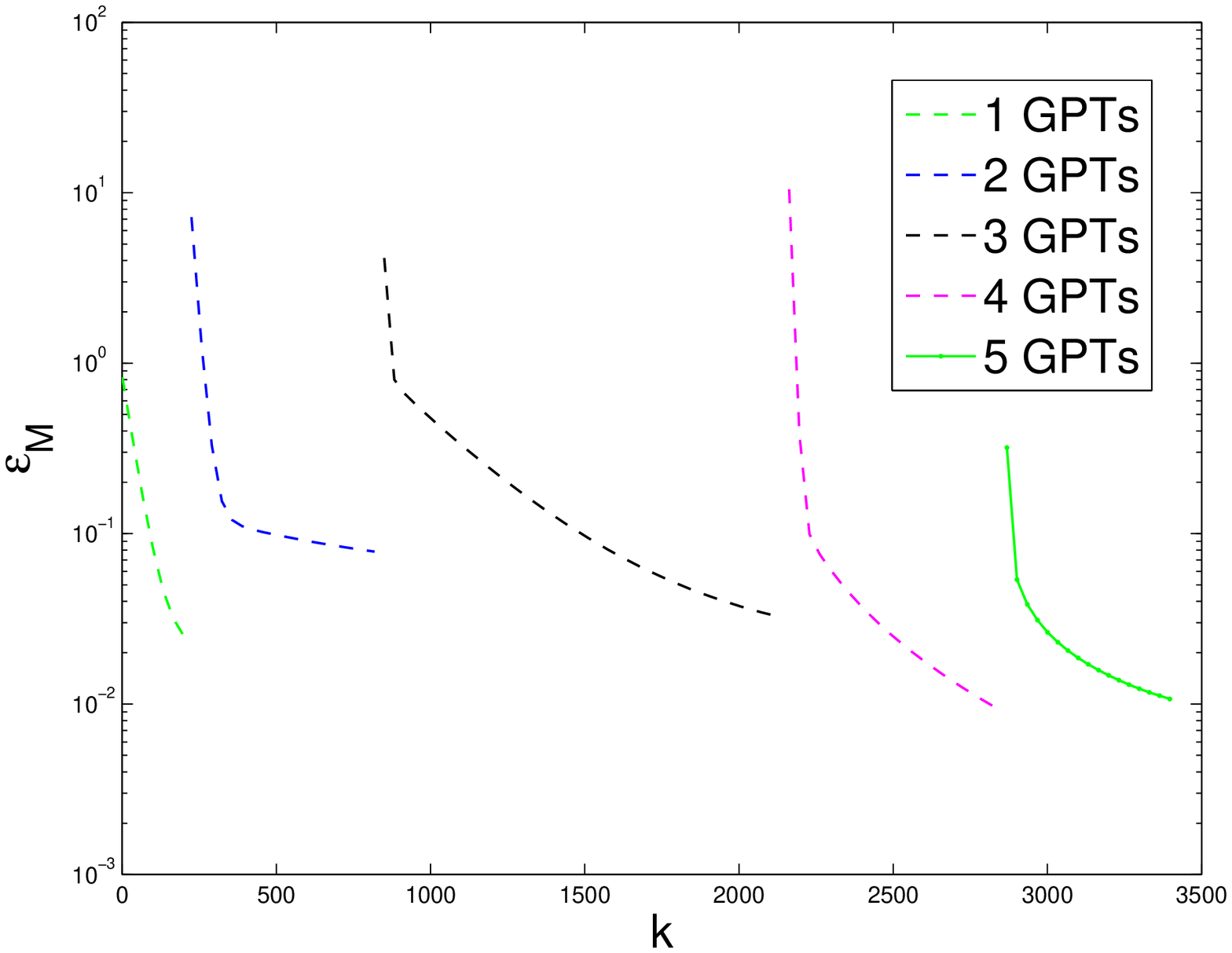}
  \includegraphics[width=2.7in]{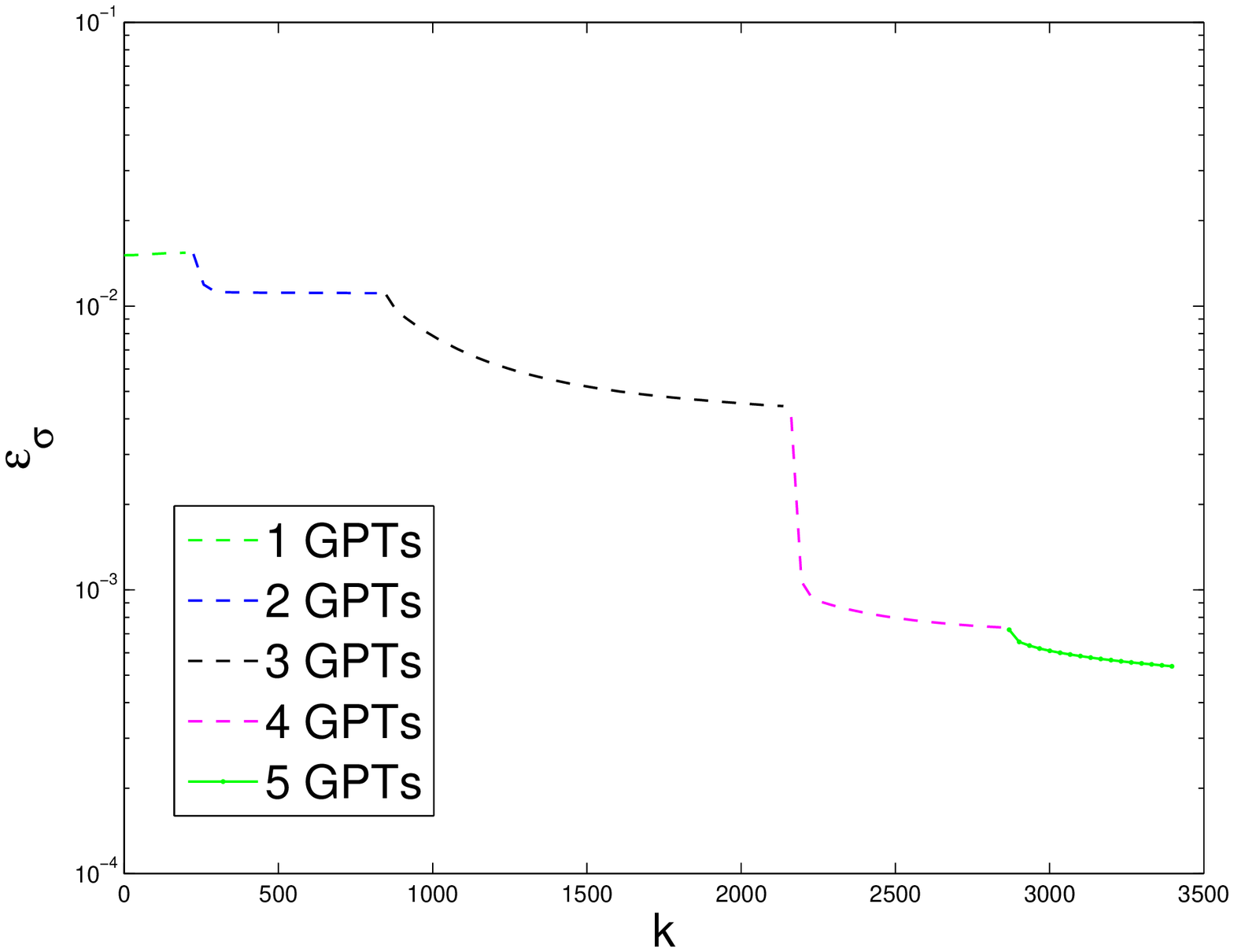}
  \end{center}
  \caption{The convergence history of $\varepsilon_{M}$ and $\varepsilon_{\Gs}$, where ${\rm k}$ is the number of iterations. $\Gs=x^3+y^5+(y-0.5)^{0.4}+3.0$. \label{fig3-result1}}
\end{figure}

\begin{figure}[ht]
\begin{center}
  \includegraphics[width=2.7in]{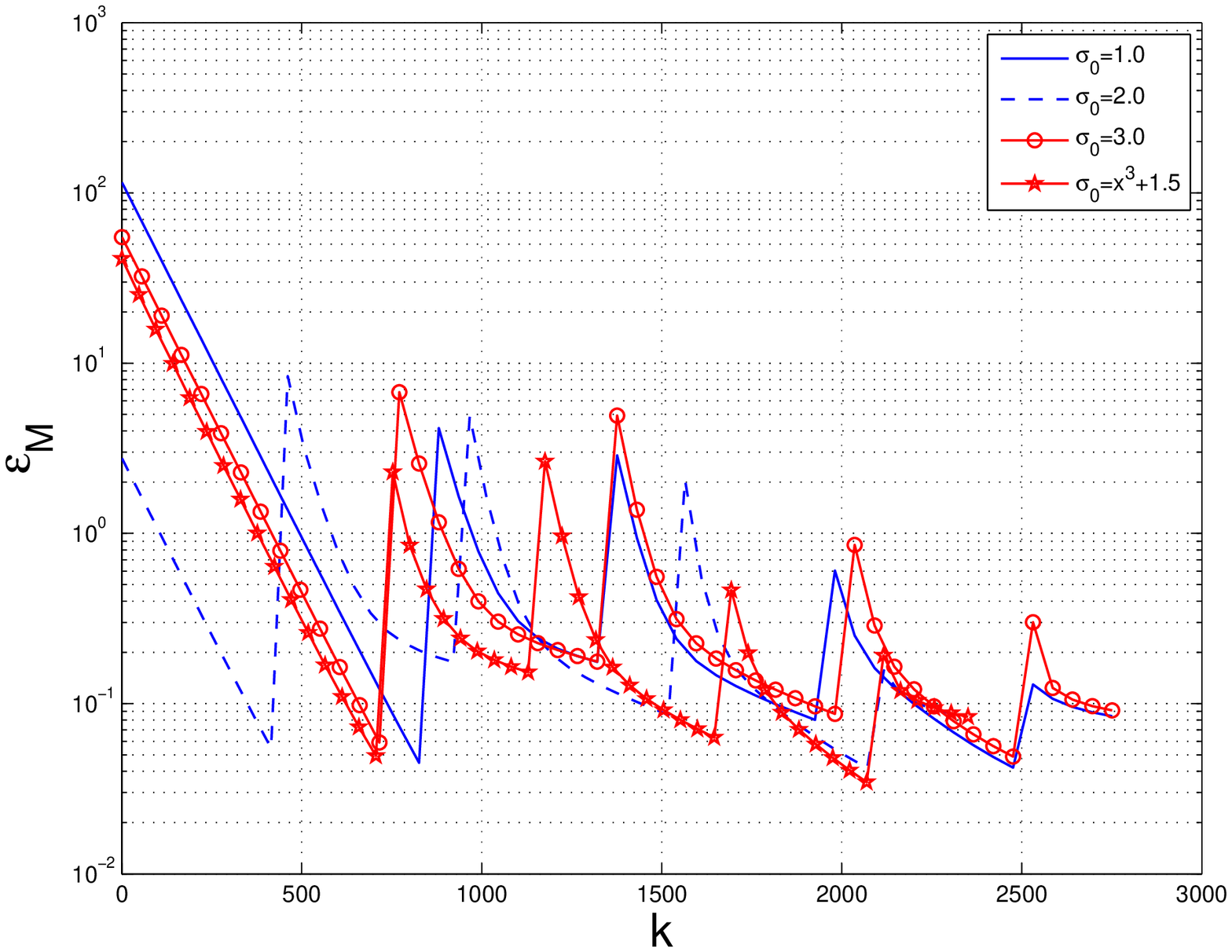}
  \includegraphics[width=2.7in]{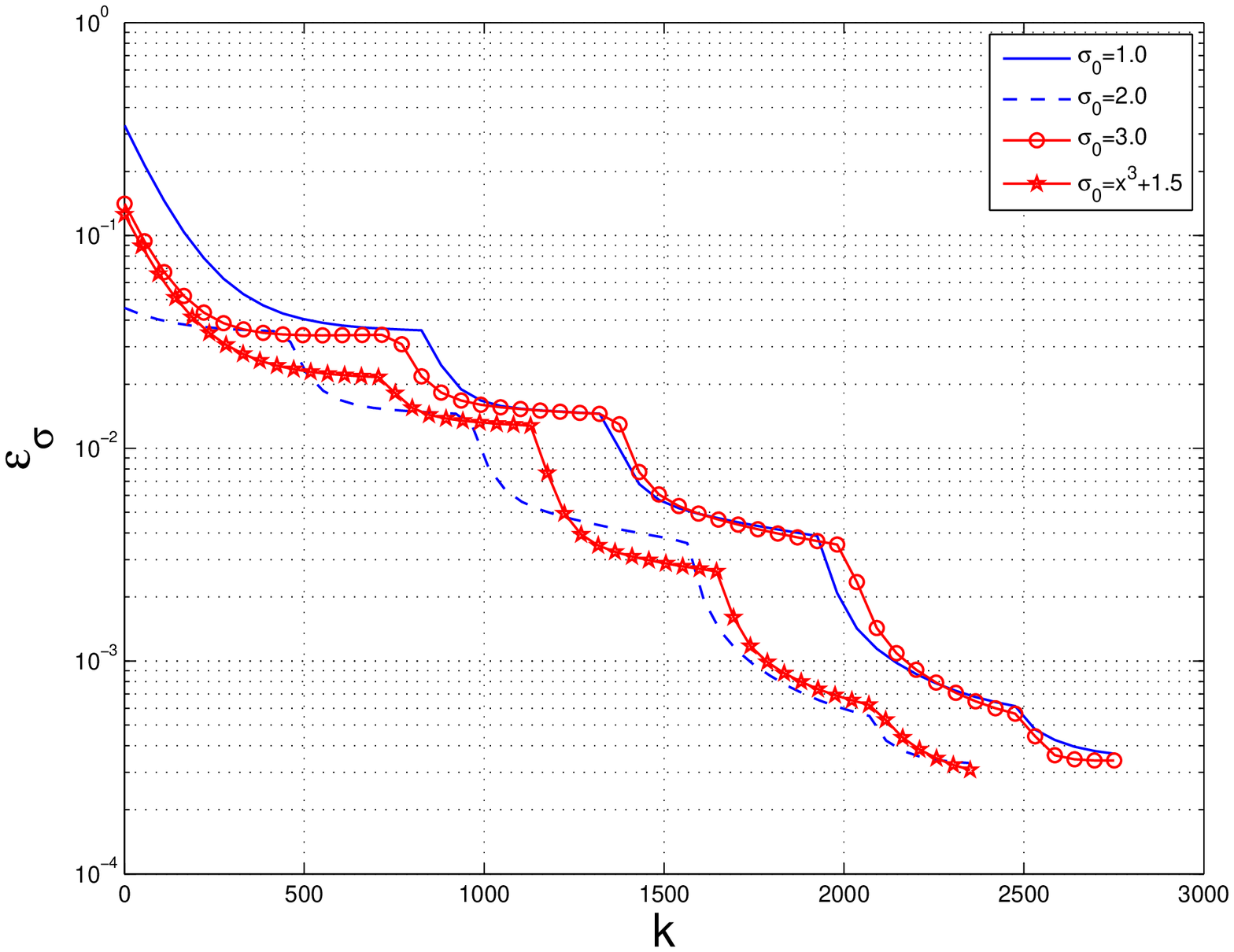}
  \end{center}
  \caption{Comparison of performance of iteration schemes by using different initial conductivities. The convergence history of $\varepsilon_{M}$ and $\varepsilon_{\Gs}$, where ${\rm k}$ is the number of iterations. $\Gs=x^3+y^5+y^2+2.0$. \label{fig4-result1}}
\end{figure}

In Table \ref{tab:1} we show the performance of the reconstruction of $C^{\infty}$ conductivity with different values of $a$ and $b$. Although
we only use quite little information from the operator $M_\Gs$ (lower order GPTs), the stability result is in some sense in accordance with
Theorem \ref{th:stability}. To test the performance of reconstruction under different initial approximation conductivity, we still use the
$C^{\infty}$ conductivity as the exact solution ($a=1.0$, $b=1.0$). Four different initial approximations are used and the reconstruction
results are shown in Figure \ref{fig4-result1}.

\begin{figure}[ht]
\begin{center}
  \includegraphics[width=1.8in]{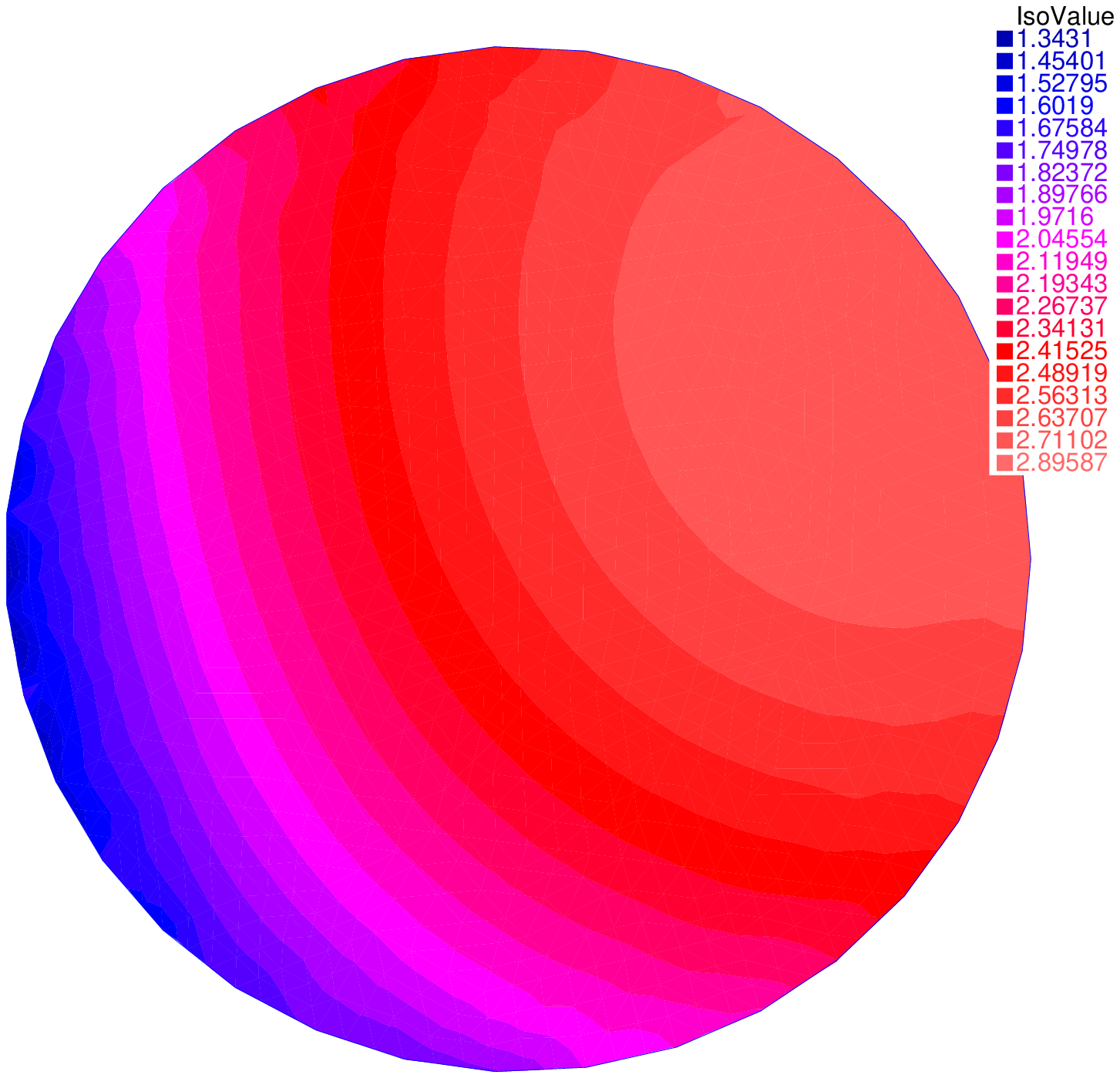}
  \includegraphics[width=1.8in]{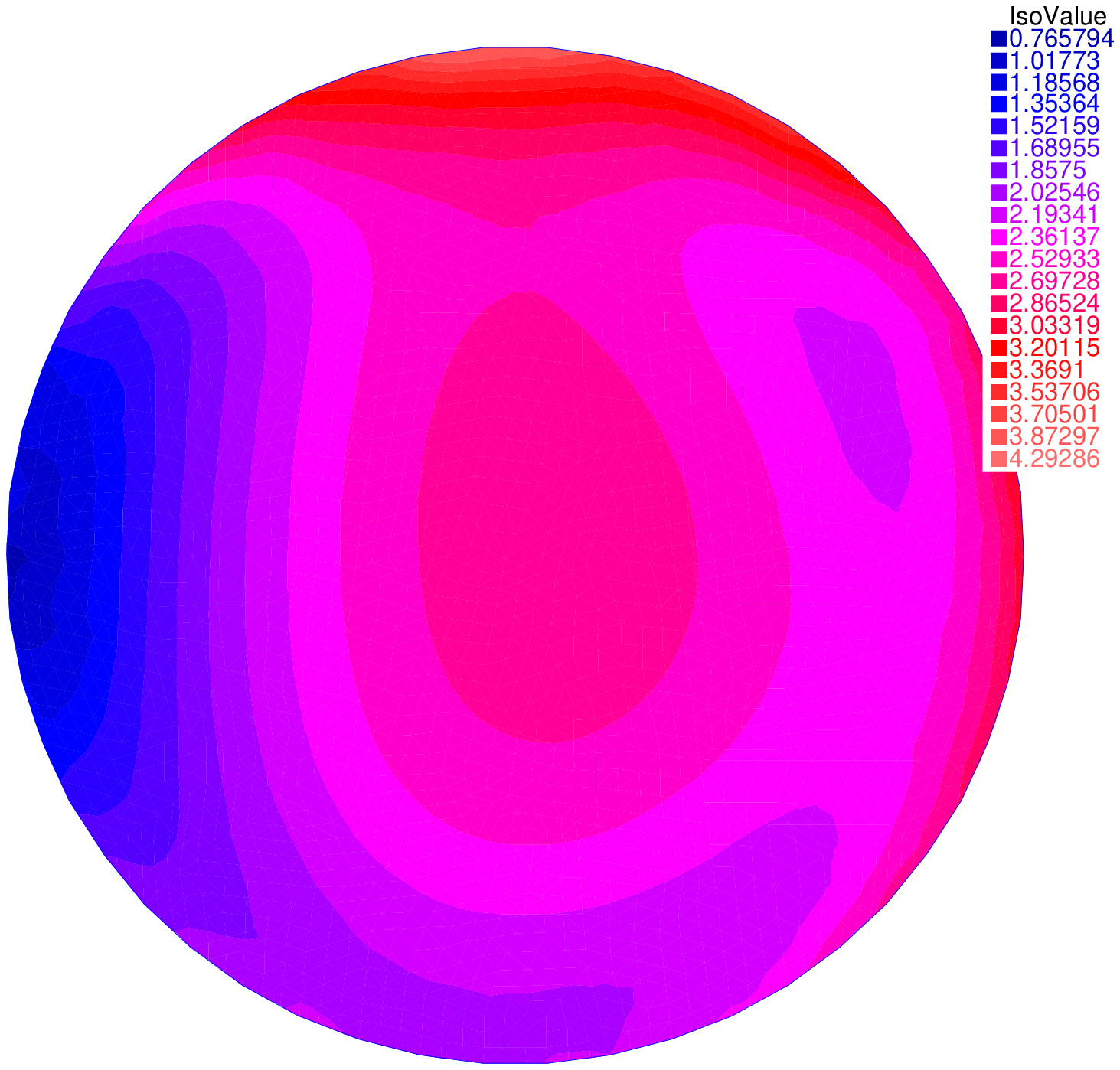}
  \includegraphics[width=1.8in]{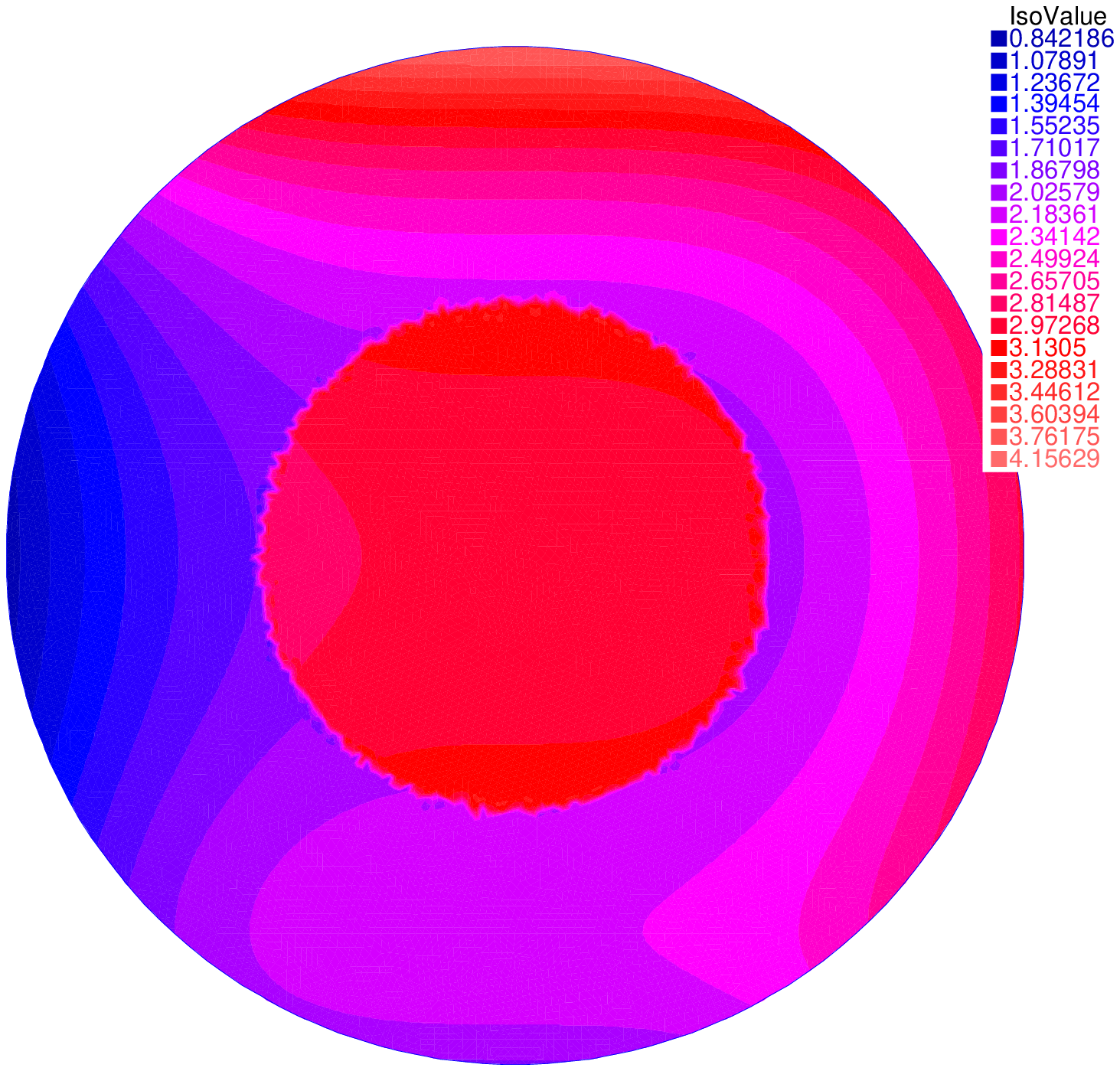}
  \end{center}
  \caption{Reconstruction of discontinuous conductivity distribution ($\Gs=x^3+y^5+y^2+\chi(x^2+y^2<0.25)+2.0$).  The first
  picture is the reconstruction result by using up to two orders of GPTs while the second one with up to five orders of
  GPTs. The third figure is the exact conductivity.\label{fig4-result}}
\end{figure}

In Figure \ref{fig4-result}, reconstruction of discontinuous conductivity (4th conductivity) by using GPTs with up to five order
is presented. It is clear that the discontinuous part is not able to be reconstructed accurately. Reconstruction of non-smooth
conductivities is quite a challenge problem, partly may be that there is no stability result for reconstruction of such
kind of conductivity. In order to uniquely determine the conductivity, higher order GPTs which contains the high
frequency information of the conductivity, must be used. As is known that higher order GPTs are quite unstable and require
a large amount of computation, which makes the reconstruction of discontinuous conductivities quite difficult. At last, we
mention that the \emph{a priori} information, or the initial approximation is especially important in reconstructing the
discontinuous conductivity. Judging from the Landweber iteration scheme, it is clear that if the initial conductivity is
smooth then all iterates shall be smooth with any order of GPTs.

\begin{table}[h]
\caption{Stability on reconstruction of different $C^{\infty}$ conductivities.}
\begin{center}
\begin{tabular}{cccccc}
\hline\noalign{\smallskip}
 $a$ & $b$ & $\mbox{iterative steps}$ & $\|\bm{Y}-\bm{M}(\Gs_{k_*})\|_F$ & $\|\Gs^*-\Gs_{k_*}\|$ & $\log\|\bm{Y}-\bm{M}(\Gs_{k_*})\|_F$ \\
 \noalign{\smallskip} \hline \noalign{\smallskip}
 0.1 & 1.0 &  1996 &  0.281857  & 0.0666685 &  -1.26635\\
 0.2 & 0.5 &  1750 &  0.265850  & 0.0604358 &  -1.32482\\
 0.2 & 1.0 &  2013 &  0.283147  & 0.0665785 &  -1.26179\\
 0.5 & 0.1 &  1802 &  0.262031  & 0.0594092 &  -1.33929\\
 0.5 & 0.5 &  1882 &  0.265512  & 0.0617185 &  -1.32610\\
 0.5 & 1.0 &  2103 &  0.285956  & 0.0676716 &  -1.25192\\
 1.0 & 0.1 &  2092 &  0.260335  & 0.0680854 &  -1.34579\\
 1.0 & 0.2 &  2100 &  0.260724  & 0.0681834 &  -1.34429\\
 1.0 & 0.5 &  2159 &  0.264578  & 0.0690950 &  -1.32962\\
 1.0 & 1.0 &  2348 &  0.289591  & 0.0738315 &  -1.23929\\
\noalign{\smallskip} \hline \noalign{\smallskip}
\end{tabular}
\end{center}
\label{tab:1}
\end{table}

\section{Conclusion}
We have presented the stability analysis for reconstruction of the inhomogeneous conductivity distribution.
We have shown the linear resolution by using optimization methods in reconstructing the conductivity.
Regularization methods are also introduced in solving the corresponding inverse problem. We have shown
how to reconstruct the GPTs by using MSR matrix in real world applications. The numerical implementations
show that the conductivities with higher smoothness and better symmetry is easier to reconstruct.
Reconstruction of discontinuous conductivities by using GPTs is much more difficult than reconstructing
the smooth conductivities and shall be considered as further works.

\end{document}